\documentclass[11pt]{article}
\usepackage[]{amsmath,amssymb}

\newtheorem{theorem}{Theorem}[section]
\newtheorem{lemma}[theorem]{Lemma}

\newtheorem{corollary}[theorem]{Corollary}
\newtheorem{exAux}[theorem]{Example}
\newenvironment{example}{\begin{exAux} \rm}{\end{exAux}}
\newtheorem{Def}[theorem]{Definition}
\newenvironment{definition}{\begin{Def} \rm}{\end{Def}}
\newtheorem{Note}[theorem]{Note}
\newenvironment{note}{\begin{Note} \rm}{\end{Note}}
\newtheorem{Problem}[theorem]{Problem}
\newenvironment{problem}{\begin{Problem} \rm}{\end{Problem}}
\newtheorem{Rem}[theorem]{Remark}

\newtheorem{Not}[theorem]{Notation}
\newenvironment{notation}{\begin{Not} \rm}{\end{Not}}
\newtheorem{Ass}[theorem]{Assumption}

\newenvironment{proof}{\medskip\noindent{\bf Proof.\ }}{\qed\medskip}

\newcommand{\qed}{\hfill\mbox{$\Box$\qquad\qquad}}

\newcommand{\Mat}[1]{\text{\rm Mat}_{#1}(\mathbb{K})}

\begin{document}
\thispagestyle{empty}

\begin{center}
\LARGE \bf
\noindent
The switching element for a Leonard pair
\end{center}

\smallskip

\begin{center}
\Large
Kazumasa Nomura and Paul Terwilliger
\end{center}

\smallskip

\begin{quote}
\small 
\begin{center}
\bf Abstract
\end{center}
Let $V$ denote a vector space with finite positive dimension.
We consider a pair of linear transformations 
$A : V \to V$ and $A^* : V \to V$ that satisfy (i) and (ii) below:
\begin{itemize}
\item[(i)] There exists a basis for $V$ with respect to which the
matrix representing $A$ is irreducible tridiagonal and the matrix
representing $A^*$ is diagonal.
\item[(ii)] There exists a basis for $V$ with respect to which the
matrix representing $A^*$ is irreducible tridiagonal and the matrix
representing $A$ is diagonal.
\end{itemize}
We call such a pair a {\em Leonard pair} on $V$.
Let $\{v_i\}_{i=0}^d$ (resp. $\{w_i\}_{i=0}^d$)
denote a basis for $V$ referred 
to in (i) (resp. (ii)).
We show that there exists a unique linear transformation
$S: V \to V$ that
sends $v_0$ to a scalar multiple of $v_d$,
fixes $w_0$, and sends $w_i$ to a 
scalar multiple of $w_i$ for $1 \leq i \leq d$.
We call $S$ the {\it switching element}.
We describe $S$ from many points of view.
\end{quote}

\section{Leonard pairs}

We begin by recalling the notion of a Leonard pair.
We will use the following terms.
A square matrix $X$ is said to be {\em tridiagonal}
whenever each nonzero entry lies on either the diagonal, the subdiagonal,
or the superdiagonal. Assume $X$ is tridiagonal.
Then $X$ is said to be {\em irreducible}
whenever each entry on the subdiagonal is nonzero and each entry on
the superdiagonal is nonzero.
We now define a Leonard pair.
For the rest of this paper $\mathbb{K}$ will denote a field.

\medskip

\begin{definition}  \cite{T:Leonard}  \label{def:LP}      \samepage
Let $V$ denote a vector space over $\mathbb{K}$ with finite positive
dimension.
By a {\em Leonard pair} on $V$ we mean an ordered pair $A,A^*$
where $A:V \to V$ and $A^*:V \to V$ are linear transformations
that satisfy (i) and (ii) below:
\begin{itemize}
\item[(i)] There exists a basis for $V$ with respect to which the
matrix representing $A$ is irreducible tridiagonal and the matrix
representing $A^*$ is diagonal.
\item[(ii)] There exists a basis for $V$ with respect to which the
matrix representing $A^*$ is irreducible tridiagonal and the matrix
representing $A$ is diagonal.
\end{itemize}
\end{definition}

\medskip

\begin{note}
It is a common notational convention to use $A^*$ to represent the
conjugate-transpose of $A$. We are {\em not} using this convention.
In a Leonard pair $A,A^*$ the linear transformations $A$ and
$A^*$ are arbitrary subject to (i) and (ii) above.
\end{note}

\medskip

We refer the reader to 
\cite{H},
\cite{N:aw},
\cite{NT:balanced}, \cite{NT:formula}, \cite{NT:det}, \cite{NT:mu},
\cite{NT:span},
\cite{P}, \cite{T:sub1}, \cite{T:sub3}, \cite{T:Leonard},
\cite{T:24points}, \cite{T:canform}, \cite{T:intro},
\cite{T:intro2}, \cite{T:split}, \cite{T:array}, \cite{T:qRacah},
\cite{T:survey}, \cite{TV}, \cite{V}, \cite{V:AW}
for background on Leonard pairs.
We especially recommend the survey \cite{T:survey}.
See \cite{AC}, \cite{AC2}, \cite{BT:Borel}, \cite{BT:loop},
\cite{H:tetra}, \cite{HT:tetra},
\cite{ITT}, \cite{IT:shape},
\cite{IT:uqsl2hat}, \cite{IT:non-nilpotent}, 
\cite{IT:tetra}, \cite{IT:inverting},
\cite{ITW:equitable}, 
\cite{N:refine}, \cite{N:height1},
\cite{R:multi}, \cite{R:6j},
\cite{T:qSerre}, \cite{T:Kac-Moody},
\cite{Z}
for related topics.

\medskip

We now give an informal summary of the present paper; our
formal treatment will begin in Section 2.
Let $A,A^*$ denote the Leonard pair on $V$ from Definition
 \ref{def:LP}.
Let $\{v_i\}_{i=0}^d$ (resp. $\{w_i\}_{i=0}^d$) denote a basis for
$V$ referred 
to in part (i) (resp. (ii)) of that definition.
We show that there exists a unique linear transformation
$S: V \to V$
that sends $v_0$ to a scalar multiple of $v_d$, 
fixes $w_0$, and  sends $w_i$ to a 
scalar multiple of $w_i$ for $1 \leq i \leq d$.
We call $S$ the {\it switching element} for $A,A^*$.
We show $S$ is invertible.
There is a well-known correspondence between Leonard pairs
and sequences of orthogonal polynomials from the terminating
branch of the Askey scheme
\cite{KoeSwa},
\cite{T:survey}; we show that $S=u_d(A)$ where 
$\{u_i\}_{i=0}^d$ are the polynomials that correspond to $A,A^*$.
A {\it flag} on $V$ is a sequence $\lbrace F_i\rbrace_{i=0}^d$
of subspaces of $V$ such that $F_i$ has dimension $i+1$
for $0 \leq i \leq d$ and $F_{i-1}\subseteq F_i$ for
$1 \leq i \leq d$. Following 
\cite[Definition 7.2]{T:24points} 
we define four flags on $V$ called 
$\lbrack 0 \rbrack$,
$\lbrack D \rbrack$,
$\lbrack 0^* \rbrack$,
$\lbrack D^* \rbrack$;
for $0 \leq i \leq d$ the 
 $i^\text{th}$ component of
$\lbrack 0 \rbrack$
(resp. $\lbrack D \rbrack$,
$\lbrack 0^* \rbrack$,
$\lbrack D^* \rbrack$) is
$\text{Span}\{w_0, w_1, \ldots, w_i\}$
(resp. 
$\text{Span}\{w_d, w_{d-1}, \ldots, w_{d-i}\}$,
$\text{Span}\{v_0, v_1, \ldots, v_i\}$,
$\text{Span}\{v_d, v_{d-1}, \ldots, v_{d-i}\}$).
These four flags are mutually opposite in the
sense of
\cite[Theorem 7.3]{T:24points}.
We show that up to multiplication by a nonzero scalar,
$S$ is the unique linear transformation on $V$ that
fixes each of 
$\lbrack 0 \rbrack$,
$\lbrack D \rbrack$ and sends
$\lbrack 0^* \rbrack$ to 
$\lbrack D^* \rbrack$.
A {\it decomposition} of $V$ is a sequence of one-dimensional
subspaces whose direct sum is $V$.
Let $x,y$ denote an ordered pair of distinct 
elements from the
set $\{ 0,D,0^*,D^*\}$.
By \cite[Theorem 8.3]{T:24points}
there exists a decomposition $\lbrack xy\rbrack$ of $V$
such that for $0 \leq i \leq d$ the
 $i^\text{th}$ component of 
$\lbrack xy \rbrack$
 is the intersection of the 
 $i^\text{th}$ component of $\lbrack x \rbrack$
and the 
 $(d-i)^\text{th}$ component of $\lbrack y \rbrack$.
We show that up to multiplication by a nonzero scalar,
$S$ is the unique linear transformation on $V$ that
sends $\lbrack 0^*0 \rbrack$ to 
 $\lbrack D^*0 \rbrack$ and
$\lbrack 0^*D \rbrack$ to 
 $\lbrack D^*D \rbrack$.
By our earlier remarks there exists a unique linear
transformation $S^*: V \to V$
that sends $w_0$ to a scalar multiple of $w_d$,
fixes $v_0$, and sends $v_i$ to a 
scalar multiple of $v_i$ for $1 \leq i \leq d$.
We show that each component of 
 $\lbrack 0^*D \rbrack$ 
(resp. 
 $\lbrack D^*D \rbrack$,
$\lbrack 0^*0 \rbrack$,
 $\lbrack D^*0 \rbrack$)
is an eigenspace for
$S^*S^{-1}S^{*-1}S$
(resp. $S^*SS^{*-1}S^{-1}$,
$S^{*-1}S^{-1}S^*S$,
$S^{*-1}SS^*S^{-1}$).
We find the corresponding eigenvalues.
We consider a certain basis for $V$ whose
 $i^\text{th}$ component is contained in
the 
 $i^\text{th}$ component of $\lbrack 0^*D\rbrack$
for $0 \leq i \leq d$. With respect to this basis
the matrix representing $A$ (resp. $A^*$)  is lower bidiagonal
(resp. upper bidiagonal) 
\cite[Lemma 3.9]{T:Leonard}.
We display the matrices
that represent $S$ and $S^*$ with respect to this basis.
In a related result we characterize the Leonard pair concept
in terms of the switching element.
We finish the paper with some open problems.

\section{Leonard systems}

When working with a Leonard pair, it is convenient to consider a closely
related object called a {\em Leonard system}. 
To prepare for our definition
of a Leonard system, we recall a few concepts from linear algebra.
Let $d$ denote a nonnegative integer and let
$\Mat{d+1}$ denote the $\mathbb{K}$-algebra consisting of all $d+1$ by
$d+1$ matrices that have entries in $\mathbb{K}$. 
We index the rows and  columns by $0, 1, \ldots, d$. 
We let $\mathbb{K}^{d+1}$ denote the $\mathbb{K}$-vector space of all
$d+1$ by $1$ matrices that have entries in $\mathbb{K}$. We index the
rows by $0,1,\ldots,d$. We view $\mathbb{K}^{d+1}$ as a left module
for $\Mat{d+1}$. We observe this module is irreducible. 
For the rest of this paper, let $\cal A$ denote a $\mathbb{K}$-algebra 
isomorphic to $\Mat{d+1}$
and let $V$ denote a simple left $\cal A$-module. We remark that $V$ is unique
up to isomorphism of $\cal A$-modules, and that $V$ has dimension $d+1$.
Let $\{v_i\}_{i=0}^d$ denote a basis for $V$.
For $X \in {\cal A}$ and $Y \in \Mat{d+1}$, we say 
{\em $Y$ represents $X$ with respect to} $\{v_i\}_{i=0}^d$
whenever $X v_j = \sum_{i=0}^d Y_{ij}v_i$ for $0 \leq j \leq d$.
For $A \in \cal A$ we say $A$ is {\em multiplicity-free}
whenever it has $d+1$ mutually distinct eigenvalues in $\mathbb{K}$. 
Assume $A$ is multiplicity-free. 
Let $\{\theta_i\}_{i=0}^d$ denote an ordering 
of the eigenvalues of $A$, and for $0 \leq i \leq d$ put
\begin{equation}        \label{eq:defEi}
    E_i = \prod_{\stackrel{0 \leq j \leq d}{j\neq i}}
             \frac{A-\theta_j I}{\theta_i - \theta_j},
\end{equation}
where $I$ denotes the identity of $\cal A$. 
We observe
(i) $AE_i = \theta_i E_i$ $(0 \leq i \leq d)$;
(ii) $E_i E_j = \delta_{i,j} E_i$ $(0 \leq i,j \leq d)$;
(iii) $\sum_{i=0}^{d} E_i = I$;
(iv) $A = \sum_{i=0}^{d} \theta_i E_i$.
Let $\cal D$ denote the subalgebra of $\cal A$ generated by $A$.
Using (i)--(iv) we find the sequence $\{E_i\}_{i=0}^d$
is a basis for the $\mathbb{K}$-vector space $\cal D$.
We call $E_i$ the {\em primitive idempotent} of $A$ associated with
$\theta_i$. 
It is helpful to think of these primitive idempotents as follows.
Observe 
\[
 V=E_0V+E_1V+\cdots+E_dV  \qquad \text{(direct sum)}.
\]
For $0 \leq i \leq d$, $E_iV$ is the (one-dimensional) eigenspace of $A$
in $V$ associated with the eigenvalue $\theta_i$, and $E_i$ acts on $V$
as the projection onto this eigenspace.
We note that for $X \in {\cal A}$ the following are equivalent:
(i) $X \in {\cal D}$; 
(ii) $XA=AX$; 
(iii) $XE_iV \subseteq E_iV$ for $0 \leq i \leq d$.

\medskip

By a {\em Leonard pair in $\cal A$} we mean an ordered pair of elements
taken from $\cal A$ that act on $V$ as a Leonard pair in the sense of
Definition \ref{def:LP}.
We now define a Leonard system.

\medskip

\begin{definition}  \cite{T:Leonard}     \label{def:LS}   \samepage
By a {\em Leonard system} in $\cal A$ we mean a sequence
\[
  \Phi= (A; \{E_i\}_{i=0}^d; A^*; \{E^*_i\}_{i=0}^d)
\]
that satisfies (i)--(v) below.
\begin{itemize}
\item[(i)] Each of $A$, $A^*$ is a multiplicity-free element in $\cal A$.
\item[(ii)] $\{E_i\}_{i=0}^d$  is an ordering of the
   primitive idempotents of $A$.
\item[(iii)] $\{E^*_i\}_{i=0}^d$ is an ordering of the
   primitive idempotents of $A^*$.
\item[(iv)] For $0 \leq i,j \leq d$, 
\begin{equation}           \label{eq:Astrid}
   E_i A^* E_j =
    \begin{cases}  
        0 & \text{\rm if $|i-j|>1$},  \\
        \neq 0 & \text{\rm if $|i-j|=1$}.
    \end{cases}
\end{equation}
\item[(v)] For $0 \leq i,j \leq d$, 
\begin{equation}             \label{eq:Atrid}
   E^*_i A E^*_j =
    \begin{cases}  
        0 & \text{\rm if $|i-j|>1$},  \\
        \neq 0 & \text{\rm if $|i-j|=1$}.
    \end{cases}
\end{equation}
\end{itemize}
We say {\em $\Phi$ is over $\mathbb{K}$}.
\end{definition}

\medskip

Leonard systems are related to Leonard pairs as follows.
Let $(A; \{E_i\}_{i=0}^d;A^*; \{E^*_i\}_{i=0}^d)$ denote a Leonard system
in $\cal A$. Then $A,A^*$ is a Leonard pair in $\cal A$
\cite[Section 3]{T:qRacah}.
Conversely, suppose $A,A^*$ is a Leonard pair in $\cal A$.
Then each of $A,A^*$ is multiplicity-free \cite[Lemma 1.3]{T:Leonard}.
Moreover there exists an ordering $\{E_i\}_{i=0}^d$ of the
primitive idempotents of $A$, and 
there exists an ordering $\{E^*_i\}_{i=0}^d$  of the
primitive idempotents of $A^*$, such that
$(A; \{E_i\}_{i=0}^d; A^*; \{E^*_i\}_{i=0}^d)$
is a Leonard system in $\cal A$ \cite[Lemma 3.3]{T:qRacah}.

\section{The $D_4$ action}

For a given Leonard system
$\Phi=(A; \{E_i\}_{i=0}^d; A^*; \{E^*_i\}_{i=0}^d)$
in $\cal A$,
each of the following is a Leonard system in $\cal A$:
\begin{eqnarray*}
\Phi^{*}  &:=& 
       (A^*; \{E^*_i\}_{i=0}^d; A; \{E_i\}_{i=0}^d),  \\
\Phi^{\downarrow} &:=&
       (A; \{E_i\}_{i=0}^d; A^*; \{E^*_{d-i}\}_{i=0}^d), \\
\Phi^{\Downarrow} &:=&
       (A; \{E_{d-i}\}_{i=0}^d; A^*; \{E^*_{i}\}_{i=0}^d).
\end{eqnarray*}
Viewing $*$, $\downarrow$, $\Downarrow$ as permutations on the set of
all Leonard systems in $\cal A$,
\begin{equation}    \label{eq:relation1}
*^2 = \downarrow^2 = \Downarrow^2 = 1,
\end{equation}
\begin{equation}    \label{eq:relation2}
\Downarrow * = * \downarrow, \quad
\downarrow * = * \Downarrow, \quad
\downarrow \Downarrow = \Downarrow \downarrow.
\end{equation}
The group generated by symbols $*$, $\downarrow$, $\Downarrow$ subject
to the relations (\ref{eq:relation1}), (\ref{eq:relation2}) is the
dihedral group $D_4$. We recall that $D_4$ is the group of symmetries of a
square and has $8$ elements.
Apparently $*$, $\downarrow$, $\Downarrow$ induce an action of $D_4$
on the set of all Leonard systems in $\cal A$.
Two Leonard systems will be called {\em relatives} whenever they are
in the same orbit of this $D_4$ action. 
The relatives of $\Phi$ are as follows:

\medskip
\noindent
\begin{center}
\begin{tabular}{c|c}
name  &  relative \\
\hline
$\Phi$ & 
       $(A; \{E_i\}_{i=0}^d; A^*;  \{E^*_i\}_{i=0}^d)$ \\ 
$\Phi^{\downarrow}$ &
       $(A; \{E_i\}_{i=0}^d; A^*;  \{E^*_{d-i}\}_{i=0}^d)$ \\ 
$\Phi^{\Downarrow}$ &
       $(A; \{E_{d-i}\}_{i=0}^d; A^*;  \{E^*_i\}_{i=0}^d)$ \\ 
$\Phi^{\downarrow \Downarrow}$ &
       $(A; \{E_{d-i}\}_{i=0}^d; A^*;  \{E^*_{d-i}\}_{i=0}^d)$ \\ 
$\Phi^{*}$  & 
       $(A^*; \{E^*_i\}_{i=0}^d; A;  \{E_i\}_{i=0}^d)$ \\ 
$\Phi^{\downarrow *}$ &
       $(A^*; \{E^*_{d-i}\}_{i=0}^d; A;  \{E_i\}_{i=0}^d)$ \\ 
$\Phi^{\Downarrow *}$ &
       $(A^*; \{E^*_i\}_{i=0}^d; A;  \{E_{d-i}\}_{i=0}^d)$ \\ 
$\Phi^{\downarrow \Downarrow *}$ &
       $(A^*; \{E^*_{d-i}\}_{i=0}^d; A;  \{E_{d-i}\}_{i=0}^d)$
\end{tabular}
\end{center}
We will use the following notational convention.

\medskip

\begin{definition}
For $g \in D_4$ and for an object $f$ associated with $\Phi$ we let
$f^g$ denote the corresponding object associated with $\Phi^{g^{-1}}$.
\end{definition}

\section{The parameter array}

In this section we recall some parameters.

\medskip

\begin{definition}        \label{def:th}
Let $\Phi=(A; \{E_i\}_{i=0}^d;A^*; \{E^*_i\}_{i=0}^d)$ denote a 
Leonard system in $\cal A$.
For $0 \leq i \leq d$ we let $\theta_i$ (resp. $\theta^*_i$)
denote the eigenvalue of $A$ (resp. $A^*$) associated with
$E_i$ (resp. $E^*_i$).
We refer to $\{\theta_i\}_{i=0}^d$ (resp. $\{\theta^*_i\}_{i=0}^d$)
as the {\em eigenvalue sequence} (resp. {\em dual eigenvalue sequence})
of $\Phi$.
We observe $\{\theta_i\}_{i=0}^d$ (resp. $\{\theta^*_i\}_{i=0}^d$)
are mutually distinct
and contained in $\mathbb{K}$.
\end{definition}

\medskip

We will use the following notation.
Let $\lambda$ denote an indeterminate and let $\mathbb{K}[\lambda]$
denote the $\mathbb{K}$-algebra consisting of all polynomials in $\lambda$
that have coefficients in $\mathbb{K}$.

\medskip

\begin{definition}             \label{def:tau}
Let $\Phi=(A; \{E_i\}_{i=0}^d; A^*; \{E^*_i\}_{i=0}^d)$
denote a Leonard system in $\cal A$.
Let $\{\theta_i\}_{i=0}^d$ (resp. $\{\theta^*_i\}_{i=0}^d$)
denote the eigenvalue sequence (resp. dual eigenvalue sequence)
of $\Phi$.
For $0 \leq i \leq d$ we define polynomials 
$\tau_i$, $\eta_i$, $\tau^*_i$, $\eta^*_i$ in $\mathbb{K}[\lambda]$
as follows:
\begin{eqnarray*}
 \tau_i &=& (\lambda-\theta_0)(\lambda-\theta_1)\cdots(\lambda-\theta_{i-1}), \\
 \eta_i &=& (\lambda-\theta_d)(\lambda-\theta_{d-1})\cdots(\lambda-\theta_{d-i+1}),\\
 \tau^*_i &=& 
   (\lambda-\theta^*_0)(\lambda-\theta^*_1)\cdots(\lambda-\theta^*_{i-1}), \\
\eta^*_i &=& 
   (\lambda-\theta^*_d)(\lambda-\theta^*_{d-1})\cdots(\lambda-\theta^*_{d-i+1}).
\end{eqnarray*}
Note that each of
$\tau_i$, $\eta_i$, $\tau^*_i$, $\eta^*_i$ is monic with degree $i$ 
for $0 \leq i \leq d$.
\end{definition}

\medskip

\begin{definition}   \cite[Theorem 4.6]{NT:formula}    \label{def:splitseq}
Let $\Phi=(A; \{E_i\}_{i=0}^d;A^*; \{E^*_i\}_{i=0}^d)$ denote a 
Leonard system in $\cal A$.
Referring to Definition \ref{def:tau}, we define scalars
\begin{eqnarray}
 \varphi_i 
   &=& (\theta^*_0-\theta^*_i)
       \frac{\text{tr}(\tau_i(A)E^*_0)}
            {\text{tr}(\tau_{i-1}(A)E^*_0)}  \qquad (1 \leq i \leq d), 
                               \label{eq:defvarphi}   \\
 \phi_i
   &=& (\theta^*_0-\theta^*_i)
       \frac{\text{tr}(\eta_i(A)E^*_0)}
            {\text{tr}(\eta_{i-1}(A)E^*_0)}  \qquad (1 \leq i \leq d),
                               \label{eq:defphi}
\end{eqnarray}
where tr means trace.
We note that in (\ref{eq:defvarphi}), (\ref{eq:defphi}) the denominators 
are nonzero by \cite[Corollary 4.5]{NT:formula}.
The sequence $\{\varphi_i\}_{i=1}^d$ 
(resp. $\{\phi_i\}_{i=1}^d$) is called the 
{\em first split sequence} (resp. {\em second split sequence}) of $\Phi$.
\end{definition}

\medskip

\begin{definition}            \label{def:param}
Let $\Phi=(A; \{E_i\}_{i=0}^d; A^*; \{E^*_i\}_{i=0}^d)$
denote a Leonard system in $\cal A$.
By the {\em parameter array of $\Phi$} we mean the sequence
$(\{\theta_i\}_{i=0}^d; \{\theta^*_i\}_{i=0}^d;
        \{\varphi_i\}_{i=1}^d; \{\phi_i\}_{i=1}^d)$,
where the $\theta_i$, $\theta^*_i$ are from
Definition \ref{def:th} and the $\varphi_i$, $\phi_i$ are
from Definition \ref{def:splitseq}.
\end{definition}

\medskip

\begin{theorem}    \cite[Theorem 1.9]{T:Leonard}   \label{thm:classify}
Let $d$ denote a nonnegative integer and let
\begin{equation}            \label{eq:paramarray}
(\{\theta_i\}_{i=0}^d; \{\theta^*_i\}_{i=0}^d;
        \{\varphi_i\}_{i=1}^d; \{\phi_i\}_{i=1}^d)
\end{equation}
denote a sequence of scalars taken from $\mathbb{K}$.
Then there exists a Leonard system $\Phi$ over $\mathbb{K}$
with parameter array (\ref{eq:paramarray}) if and only if
(PA1)--(PA5) hold below.
\begin{itemize}
\item[]
\begin{itemize}
\item[(PA1)]  $\varphi_i \neq 0$, $\phi_i \neq 0$ $(1 \leq i \leq d)$.
\item[(PA2)]  $\theta_i \neq \theta_j$, $\theta^*_i \neq \theta^*_j$
   if $i \neq j$ $(0 \leq i,j \leq d)$.
\item[(PA3)] For $1 \leq i \leq d$,
\[
 \varphi_i = \phi_1 
 \sum_{h=0}^{i-1} \frac{\theta_h - \theta_{d-h}}
                       {\theta_0 - \theta_d}
         + (\theta^*_{i}-\theta^*_{0})(\theta_{i-1}-\theta_{d}).
\]
\item[(PA4)] For $1 \leq i \leq d$,
\[
\phi_i = \varphi_1
\sum_{h=0}^{i-1} \frac{\theta_h - \theta_{d-h}}
                       {\theta_0 - \theta_d}
         + (\theta^*_{i}-\theta^*_{0})(\theta_{d-i+1}-\theta_{0}).
\]
\item[(PA5)]  The expressions
\begin{equation}        \label{eq:indep}
   \frac{\theta_{i-2}-\theta_{i+1}}{\theta_{i-1}-\theta_{i}},
 \quad
   \frac{\theta^*_{i-2}-\theta^*_{i+1}}{\theta^*_{i-1}-\theta^*_{i}}
\end{equation}
are equal and independent of $i$ for $2 \leq i \leq d-1$.
\end{itemize}
\end{itemize}
Suppose (PA1)--(PA5) hold. Then $\Phi$ is unique up to isomorphism
of Leonard systems.
\end{theorem}

\medskip

The $D_4$ action affects the parameter array as follows.

\medskip

\begin{lemma}   \cite[Theorem 1.11]{T:Leonard}        \label{lem:D4}   \samepage
Let $\Phi=(A; \{E_i\}_{i=0}^d; A^*; \{E^*_i\}_{i=0}^d)$
denote a Leonard system in $\cal A$ and let
$(\{\theta_i\}_{i=0}^d; \{\theta^*_i\}_{i=0}^d;
        \{\varphi_i\}_{i=1}^d; \{\phi_i\}_{i=1}^d)$
denote the parameter array of $\Phi$.
Then the following (i)--(iii) hold.
\begin{itemize}
\item[(i)] The parameter array of $\Phi^*$ is 
\[
   (\{\theta^*_i\}_{i=0}^d; \{\theta_i\}_{i=0}^d;
        \{\varphi_i\}_{i=1}^d; \{\phi_{d-i+1}\}_{i=1}^d).
\]
\item[(ii)]
The parameter array of $\Phi^{\downarrow}$ is
\[
  (\{\theta_i\}_{i=0}^d; \{\theta^*_{d-i}\}_{i=0}^d;
        \{\phi_{d-i+1}\}_{i=1}^d; \{\varphi_{d-i+1}\}_{i=1}^d).
\]
\item[(iii)]
The parameter array of $\Phi^{\Downarrow}$ is
\[
  (\{\theta_{d-i}\}_{i=0}^d; \{\theta^*_{i}\}_{i=0}^d;
        \{\phi_{i}\}_{i=1}^d; \{\varphi_{i}\}_{i=1}^d).
\]
\end{itemize}
\end{lemma}

\medskip

We finish this section with a comment.

\medskip

\begin{lemma}           \label{lem:XrEs0Xs}
Let $(A; \{E_i\}_{i=0}^d; A^*; \{E^*_i\}_{i=0}^d)$
denote a Leonard system in $\cal A$.
Let $\cal D$ denote the subalgebra of $\cal A$ generated by $A$, and
let $X$ denote an element of $\cal D$ such that $XE^*_0=0$.
Then $X=0$.
\end{lemma}

\begin{proof}
Immediate from  \cite[Lemma 5.9]{T:qRacah}.
\end{proof}

\section{The switching element $S$}

\begin{definition}         \label{def:S}
For a Leonard system $\Phi=(A; \{E_i\}_{i=0}^d; A^*; \{E^*_i\}_{i=0}^d)$
in $\cal A$
we define 
\begin{equation}            \label{eq:defS}
  S=\sum_{r=0}^d
    \frac{\phi_d\phi_{d-1}\cdots\phi_{d-r+1}}
         {\varphi_1\varphi_2\cdots\varphi_r} E_r,
\end{equation}
where $\{\varphi_i\}_{i=1}^d$ (resp. $\{\phi_i\}_{i=1}^d$)
denotes the first (resp. second) split sequence of $\Phi$.
We call $S$ the {\em switching element} for $\Phi$.
\end{definition}

\medskip

\begin{note}     
Let $\Phi=(A; \{E_i\}_{i=0}^d; A^*; \{E^*_i\}_{i=0}^d)$
denote a Leonard system.
In what follows we will often make use of the switching element of $\Phi^*$.
By (\ref{eq:defS}) and Lemma \ref{lem:D4}(i),
\begin{equation}            \label{eq:defSs}
  S^*=\sum_{r=0}^d
    \frac{\phi_1\phi_{2}\cdots\phi_{r}}
         {\varphi_1\varphi_2\cdots\varphi_r} E^*_r.
\end{equation}
We call $S^*$  the {\em dual switching element} for $\Phi$.
\end{note}

\medskip

\begin{lemma}              \label{note:Sinv}
The switching element (\ref{eq:defS}) and the dual switching
element (\ref{eq:defSs}) are invertible with 
\begin{eqnarray}
\label{sinv}
  S^{-1} &=& \sum_{r=0}^d
    \frac{\varphi_1\varphi_2\cdots\varphi_r}
         {\phi_d\phi_{d-1}\cdots\phi_{d-r+1}} E_r,  \\
\label{ssinv}
  S^{*-1} &=& \sum_{r=0}^d
    \frac{\varphi_1\varphi_2\cdots\varphi_r}
         {\phi_1\phi_{2}\cdots\phi_{r}}   E^*_r.
\end{eqnarray}
\end{lemma}

\begin{proof} To obtain (\ref{sinv}) we note that
the sum on the right in 
(\ref{eq:defS}) times the
sum on the right in 
(\ref{sinv}) is equal to the identity; this is verified
using equations (ii), (iii) below 
 (\ref{eq:defEi}).
Line (\ref{ssinv}) is similarly obtained.
\end{proof}

\medskip

\begin{theorem}                        \label{thm:Srelative}
Let $\Phi=(A; \{E_i\}_{i=0}^d; A^*; \{E^*_i\}_{i=0}^d)$ 
denote a Leonard system  with switching element $S$ and
dual switching element $S^*$.
Then the switching element and the dual switching element
for the relatives of $\Phi$ are given in the following table:
{\small 
\[
\begin{array}{c|cccccccc}
 \text{relative} 
 & \Phi & \Phi^{\downarrow} & \Phi^{\Downarrow} & \Phi^{\downarrow\Downarrow}
 & \Phi^* & \Phi^{\downarrow*}&\Phi^{\Downarrow*}&\Phi^{\downarrow\Downarrow*}
\\
\hline
\hline
 \text{switching element}
 & S & S^{-1} & \varphi\phi^{-1} S & \varphi^{-1}\phi S^{-1}
 & S^* &  \varphi\phi^{-1} S^* & S^{*-1} & \varphi^{-1}\phi S^{*-1} 
\\
\hline
 \text{dual switching element}
 & S^* &  \varphi\phi^{-1} S^* & S^{*-1} & \varphi^{-1}\phi S^{*-1} 
 & S & S^{-1} &  \varphi\phi^{-1} S & \varphi^{-1}\phi S^{-1}
\end{array}
\]
}
In the above table we abbreviate
\[
  \varphi = \varphi_1\varphi_2\cdots\varphi_d,
  \qquad \qquad
  \phi = \phi_1\phi_2\cdots\phi_d,
\]
where $\{\varphi_i\}_{i=1}^d$ (resp. $\{\phi_i\}_{i=1}^d$) is the 
first (resp. second) split sequence of $\Phi$.
\end{theorem}

\begin{proof}
Apply $D_4$ to (\ref{eq:defS}) and use Lemma \ref{lem:D4}.
\end{proof}

\medskip

We now describe the switching element from various points of view.

\section{Representing $S$ as a polynomial}

Let $(A; \{E_i\}_{i=0}^d; A^*; \{E^*_i\}_{i=0}^d)$ denote a Leonard system 
in $\cal A$ with switching element $S$. 
Let $\cal D$ denote the subalgebra of $\cal A$ generated by
$A$, and recall $\{E_i\}_{i=0}^d$ is a basis for $\cal D$. 
By this and (\ref{eq:defS}) we find $S \in {\cal D}$, 
so $S$ is a polynomial in $A$. In the present
section we find this polynomial.

\medskip

\begin{lemma}              \label{lem:defpi}
Let $(A; \{E_i\}_{i=0}^d; A^*; \{E^*_i\}_{i=0}^d)$
denote a Leonard system in $\cal A$.
Then for $0 \leq i \leq d$ there exists a unique monic polynomial $p_i$
in $\mathbb{K}[\lambda]$ with degree $i$ such that
\[
       p_i(A)E^*_0V = E^*_iV.
\]
\end{lemma}

\begin{proof}
The existence of $p_i$ is established in \cite[Theorem 8.3]{T:qRacah}.
Concerning uniqueness, suppose we are given a monic polynomial
$p'_i$ in $\mathbb{K}[\lambda]$ of degree $i$ such that
$p'_i(A)E^*_0V = E^*_iV$.
We show $p_i=p'_i$.
To this end we define $f=p_i - p'_i$ and show $f=0$.
By construction $f(A)E^*_0V \subseteq E^*_iV$.
Each of $p_i$, $p'_i$ is monic of degree $i$ so the degree
of $f$ is at most $i-1$.
By this and (\ref{eq:Atrid}) we find $f(A)E^*_0V$ is included in
$\sum_{k=0}^{i-1} E^*_kV$.
By these comments we find $f(A)E^*_0V=0$ so $f(A)E^*_0=0$.
Now $f(A)=0$ in view of Lemma \ref{lem:XrEs0Xs}.
This implies $f=0$ since $I,A,A^2,\ldots,A^d$ are linearly independent.
\end{proof}

\medskip

\begin{lemma}       \cite[Lemma 17.5]{T:qRacah}
Let $(A; \{E_i\}_{i=0}^d; A^*; \{E^*_i\}_{i=0}^d)$
denote a Leonard system and let
$(\{\theta_i\}_{i=0}^d; \{\theta^*_i\}_{i=0}^d;
        \{\varphi_i\}_{i=1}^d; \{\phi_i\}_{i=1}^d)$
denote the corresponding parameter array.
Let the polynomials $\{p_i\}_{i=0}^d$ be from Lemma \ref{lem:defpi}.
Then 
\begin{equation}                \label{eq:pith0}
 p_i(\theta_0)=
   \frac{\varphi_1\varphi_2\cdots\varphi_i}
        {\tau^*_i(\theta^*_i)}
    \qquad (0 \leq i \leq d).
\end{equation}
Moreover $p_i(\theta_0) \neq 0$.
\end{lemma}

\medskip

\begin{definition}  \cite[Definition 14.1]{T:qRacah}  \label{def:ui}
Let $(A; \{E_i\}_{i=0}^d; A^*; \{E^*_i\}_{i=0}^d)$
denote a Leonard system over $\mathbb{K}$ and let
the polynomials $\{p_i\}_{i=0}^d$ be as in Lemma \ref{lem:defpi}.
For $0 \leq i \leq d$ we define
\begin{equation}           \label{eq:defui}
   u_i = \frac{p_i}{p_i(\theta_0)},
\end{equation}
where $\theta_0$ is from Definition \ref{def:th}.
\end{definition}

\medskip

\begin{lemma}                \label{lem:uiAEs0V}
Let $(A; \{E_i\}_{i=0}^d; A^*; \{E^*_i\}_{i=0}^d)$
denote a Leonard system in $\cal A$ and
let the polynomials $\{u_i\}_{i=0}^d$ be from Definition \ref{def:ui}.
Then 
\[
   u_i(A)E^*_0V=E^*_iV  \qquad (0 \leq i \leq d).
\]
\end{lemma}

\begin{proof}
Combine Lemma \ref{lem:defpi} and (\ref{eq:defui}).
\end{proof}

\medskip

\begin{lemma}   \cite[Theorem 14.7]{T:qRacah}    \label{lem:AWduality}
Let $(A; \{E_i\}_{i=0}^d; A^*; \{E^*_i\}_{i=0}^d)$
denote a Leonard system with eigenvalue sequence
$\{\theta_i\}_{i=0}^d$ and dual eigenvalue sequence
$\{\theta^*_i\}_{i=0}^d$.
Let the polynomials $\{u_i\}_{i=0}^d$ be as in
Definition \ref{def:ui} and recall $\{u^*_i\}_{i=0}^d$ are the
corresponding polynomials for $\Phi^*$.
Then for $0 \leq i,j \leq d$,
\begin{equation}       \label{eq:AWduality}
   u_i(\theta_j) = u^*_j(\theta^*_i).
\end{equation}
\end{lemma}

\medskip

\begin{theorem}         \label{thm:SudA}
Let $\Phi=(A; \{E_i\}_{i=0}^d; A^*; \{E^*_i\}_{i=0}^d)$
denote a Leonard system  and let
$S$ denote the switching element for $\Phi$.
Then 
\begin{equation}             \label{eq:SudA}
    S = u_d(A),
\end{equation}
where the polynomial $u_d$ is from Definition \ref{def:ui}.
\end{theorem}

\begin{proof}
By $D_4$ symmetry it suffices to show $S^*=u^*_d(A^*)$.
Using the comments below (\ref{eq:defEi}) we find
\begin{equation}        \label{eq:SudAaux1}
   u^*_d(A^*)= \sum_{i=0}^d u^*_d(\theta^*_i) E^*_i.
\end{equation}
For $0 \leq i \leq d$ we compute $u^*_d(\theta^*_i)$ as follows.
By Lemma \ref{lem:defpi} the polynomial $p_i$ is invariant under
$\Downarrow$; that is $p_i^\Downarrow = p_i$.
We apply $\Downarrow$ to (\ref{eq:pith0}) using this and
Lemma \ref{lem:D4} to get
\[
   p_i(\theta_d) =
    \frac{\phi_1\phi_2\cdots\phi_i}
         {\tau^*_i(\theta^*_i)}.
\]
Combining this with (\ref{eq:pith0}), (\ref{eq:defui}) we find
\[
   u_i(\theta_d) =
    \frac{\phi_1\phi_2\cdots\phi_i}
         {\varphi_1\varphi_2\cdots\varphi_i}.
\]
By this and Lemma \ref{lem:AWduality} we get 
\begin{equation}           \label{eq:usdthsi}
   u^*_d(\theta^*_i) =
     \frac{\phi_1\phi_2\cdots\phi_i}
          {\varphi_1\varphi_2\cdots\varphi_i}.
\end{equation}
Evaluating (\ref{eq:SudAaux1}) using (\ref{eq:usdthsi}) 
and comparing the result with (\ref{eq:defSs}) we find
$S^*=u^*_d(A^*)$.
The result follows.
\end{proof}

\medskip

The switching element is characterized as follows.

\medskip

\begin{theorem}            \label{thm:XEs0V}   \samepage
Let $(A; \{E_i\}_{i=0}^d; A^*; \{E^*_i\}_{i=0}^d)$
denote a Leonard system in $\cal A$ with switching element $S$.
Let $\cal D$ denote the subalgebra of $\cal A$
generated by $A$.
Then for all nonzero $X \in {\cal A}$
the following (i), (ii) are equivalent.
\begin{itemize}
\item[(i)] $X$ is a scalar multiple of $S$.
\item[(ii)] $X \in {\cal D}$ and $XE^*_0V \subseteq E^*_dV$.
\end{itemize}
Suppose (i), (ii) hold. Then $XE^*_0V = E^*_dV$.
\end{theorem}

\begin{proof}
(i)$\Rightarrow$(ii): 
We mentioned in the first paragraph of this section that
$S \in {\cal D}$. 
We have $SE^*_0V=E^*_dV$
by Lemma \ref{lem:uiAEs0V} and Theorem \ref{thm:SudA}.

(ii)$\Rightarrow$(i):
For $0 \leq i \leq d$ we define ${\cal D}_i = \text{Span}\{u_i(A)\}$.
Observe that
${\cal D} = \sum_{i=0}^d {\cal D}_i$ (direct sum).
Also observe by Lemma \ref{lem:uiAEs0V}
that ${\cal D}_i E^*_0V = E^*_iV$ for $0 \leq i \leq d$.
We assume $XE^*_0V \subseteq E^*_dV$ so $X \in {\cal D}_d$ and
in other words $X$ is a scalar multiple of $u_d(A)$.
By this and Theorem \ref{thm:SudA} we find $X$ is a scalar multiple
of $S$.

Now suppose (i), (ii) hold. We mentioned in the proof of
(i)$\Rightarrow$(ii) that $SE^*_0V = E^*_dV$.
But $X$ is nonzero and a scalar multiple of $S$ so
$XE^*_0V = E^*_dV$.
\end{proof}

\section{Decompositions and flags}

In this section we recall the notion of a decomposition and a flag.

\medskip

By a {\em decomposition} of $V$ we mean a sequence
$\{V_i\}_{i=0}^d$ of subspaces of $V$ such that
$V_i$ has dimension $1$ for $0 \leq i \leq d$ and
$\sum_{i=0}^d V_i=V$ (direct sum).
Let $\{V_i\}_{i=0}^d$ denote a decomposition of $V$.
By the {\em inversion} of this decomposition we mean
the decomposition $\{V_{d-i}\}_{i=0}^d$.

By a {\em flag} on $V$ we mean a sequence
$\{F_i\}_{i=0}^d$ of subspaces of $V$ such that 
$F_i$ has dimension $i+1$ for $0 \leq i \leq d$ and
$F_{i-1} \subseteq F_{i}$ for $1 \leq i \leq d$.
The following construction yields a flag on $V$.
Let $\{V_i\}_{i=0}^d$ denote a decomposition of $V$. Define
\[
   F_i = V_0+V_1+\cdots+V_i  \qquad (0 \leq i \leq d).
\]
Then $\{F_i\}_{i=0}^d$ is a flag on $V$. 
We say this flag is {\em induced} by the decomposition $\{V_i\}_{i=0}^d$.

We recall what it means for two flags on V to be {\em opposite}.
Suppose we are given two flags on $V$: 
$\{F_i\}_{i=0}^d$ and $\{F'_i\}_{i=0}^d$.
We say these flags are {\em opposite} whenever there exists a decomposition
$\{V_i\}_{i=0}^d$ of $V$ such that
\[
   F_i = V_0+V_1+\cdots+V_i, \qquad
   F'_i = V_d+V_{d-1}+\cdots+V_{d-i}
\]
for $0 \leq i \leq d$.
In this case
\[
    F_i \cap F'_j = 0 \qquad \text{if $i+j<d$}  \qquad (0 \leq i,j \leq d)
\]
and
\[
   V_i = F_i \cap F'_{d-i} \qquad (0 \leq i \leq d).
\]
In particular the decomposition $\{V_i\}_{i=0}^d$ is uniquely determined
by the given flags. We say this decomposition is
{\em induced} by the given flags.

\medskip

We end this section with some notation.

\medskip

\begin{notation}
Let $F=\{F_i\}_{i=0}^d$ denote a sequence of subspaces of $V$. 
Then for $X \in {\cal A}$ we write $XF$ to denote the sequence
$\{XF_i\}_{i=0}^d$.
We say $X$ {\em fixes} $F$ whenever $XF=F$.
Let $F'=\{F'_i\}_{i=0}^d$ denote a second sequence of subspaces of $V$.
We write $F \subseteq F'$ whenever $F_i \subseteq F'_i$ for $0 \leq i \leq d$.
\end{notation}

\section{Some decompositions and flags associated with a Leonard system}

We now return our attention to Leonard systems.
Let $\Phi=(A; \{E_i\}_{i=0}^d; A^*; \{E^*_i\}_{i=0}^d)$
denote a Leonard system in $\cal A$.
Using $\Phi$ we construct
four mutually opposite flags and consider the
decompositions that they induce.
We start with a definition.

\medskip

\begin{definition}          \label{def:Omega}
For notational convenience let
$\Omega$ denote the set consisting of four symbols
$0,D,0^*,D^*$.
\end{definition}

\medskip

\begin{definition}          \label{def:flags}
Let $(A; \{E_i\}_{i=0}^d; A^*; \{E^*_i\}_{i=0}^d)$
denote a Leonard system in $\cal A$. 
For $z \in \Omega$ we define a flag on $V$ which we denote
by $[z]$.
To define this flag we display the $i^\text{th}$ component
for $0 \leq i \leq d$.
\[
\begin{array}{c|c}
    z & \text{$i^\text{th}$ component of $[z]$}  \\
  \hline
    0 & E_0V+E_1V+\cdots+E_iV  \\
    D & E_dV+E_{d-1}V+\cdots+E_{d-i}V \\
   0^* & E^*_0V+E^*_1V+\cdots+E^*_iV \\
   D^* & E^*_dV+ E^*_{d-1}V+\cdots+E^*_{d-i}V
\end{array}
\]
\end{definition}

\medskip

\begin{lemma}               \label{lem:flagcomponent}   \samepage
Referring to Definition \ref{def:flags},
the following (i)--(iv) hold for $0 \leq i \leq d$.
\begin{itemize}
\item[(i)]  The $i^\text{th}$ component of $[0]$ is equal to $\eta_{d-i}(A)V$.
\item[(ii)] The $i^\text{th}$ component of $[D]$ is equal to $\tau_{d-i}(A)V$.
\item[(iii)] The $i^\text{th}$ component of $[0^*]$ is equal to 
$\eta^*_{d-i}(A^*)V$.
\item[(iv)] The $i^\text{th}$ component of $[D^*]$ is equal to 
$\tau^*_{d-i}(A^*)V$.
\end{itemize}
\end{lemma}

\begin{proof}
(i): Recall that $V=\sum_{j=0}^d E_jV$ (direct sum).
Further recall that for $0 \leq j \leq d$, $E_jV$ is
an eigenspace for $A$ with eigenvalue $\theta_j$.
This implies that for $0 \leq j,k \leq d$,
$(A-\theta_kI)E_jV$ equals $0$ if $j=k$ and $E_jV$ if $j \neq k$.
By these comments and Definition 4.2
we have $\eta_{d-i}(A)E_jV=0$ for $i+1 \leq j \leq d$
and $\eta_{d-i}(A)E_jV=E_jV$ for $0 \leq j \leq i$.
Therefore $\eta_{d-i}(A)V=\sum_{j=0}^i E_jV$ and this is the $i^\text{th}$
component of $[0]$.

(ii)--(iv): Similar.
\end{proof}

\medskip

\begin{lemma}  \cite[Theorem 7.3]{T:24points}     \label{lem:opposite}
The four flags in Definition \ref{def:flags} are mutually opposite.
\end{lemma}

\medskip

\begin{definition}          \label{def:decompositions}
Let $\Phi=(A; \{E_i\}_{i=0}^d; A^*; \{E^*_i\}_{i=0}^d)$
denote a Leonard system in $\cal A$.
Let $z,w$ denote an ordered pair of distinct elements of $\Omega$.
By Lemma \ref{lem:opposite} the flags $[z]$, $[w]$ are opposite.
Let $[zw]$ denote the decomposition of $V$ induced by
$[z]$, $[w]$.
\end{definition}

\medskip

We mention a few basic properties
of the decompositions from Definition \ref{def:decompositions}.

\medskip

\begin{lemma}                  \label{lem:induce}
Referring to Definition \ref{def:decompositions},
for distinct $z$, $w \in \Omega$ the following (i)--(iii) hold.
\begin{itemize} 
\item[(i)] The decomposition $[zw]$ is the inversion of $[wz]$.
\item[(ii)] For $0 \leq i \leq d$
the $i^\text{th}$ component of $[zw]$ is the intersection
of the $i^\text{th}$  component of $[z]$ and the $(d-i)^\text{th}$ 
component of $[w]$.
\item[(iii)] The decomposition $[zw]$ induces $[z]$ and 
 the inversion of $[zw]$ induces $[w]$.
\end{itemize}
\end{lemma}

\begin{proof}
Routine using Section 7 and  
Definition \ref{def:decompositions}.
\end{proof}

\medskip

\begin{example}   \label{exam:decompositions}
We display some of the decompositions from 
Definition \ref{def:decompositions}.
For each decomposition in the table below we give
the $i^\text{th}$ component for $0 \leq i \leq d$.
\[
\begin{array}{c|c}
 \text{decomposition} & \text{$i^\text{th}$ component}  \\
\hline
  \;\;\; [0^*D] \;\;\; & 
   \;\;\; (E^*_0V+\cdots+E^*_iV)\cap(E_iV+\cdots+E_dV)  \\
  {[D^*D]} & 
    (E^*_dV+\cdots+E^*_{d-i}V)\cap(E_iV+\cdots+E_dV) \\
  {[0^*0]}  &
     (E^*_0V+\cdots+E^*_{i}V)\cap(E_{d-i}V+\cdots+E_{0}V) \\
  {[D^*0]}  & 
     (E^*_dV+\cdots+E^*_{d-i}V)\cap(E_{d-i}V+\cdots+E_{0}V) \\
  {[0D]}  &
       E_iV \\
  {[0^*D^*]} &
       E^*_iV 
\end{array}
\]
\end{example}

\medskip

\begin{lemma}         \label{lem:split}
Referring to Definition \ref{def:decompositions},
the following (i)--(iv) hold for $0 \leq i \leq d$.
\begin{itemize}
\item[(i)] The $i^\text{th}$ component of $[0^*D]$ is equal to
             $\tau_i(A)E^*_0V$ and $\eta^*_{d-i}(A^*)E_dV$. 
\item[(ii)] The $i^\text{th}$ component of $[D^*D]$ is equal to 
             $\tau_i(A)E^*_dV$  and $\tau^*_{d-i}(A^*)E_dV$.
\item[(iii)] The $i^\text{th}$ component of $[0^*0]$ is equal to 
             $\eta_i(A)E^*_0V$  and $\eta^*_{d-i}(A^*)E_0V$.
\item[(iv)] The $i^\text{th}$ component of $[D^*0]$ is equal to 
             $\eta_i(A)E^*_dV$  and $\tau^*_{d-i}(A^*)E_0V$.
\end{itemize}
\end{lemma}

\begin{proof}
(i): We first show that $\tau_i(A)E^*_0V$ is equal to the $i^\text{th}$
component of $[0^*D]$.
Denote this $i^\text{th}$ component by $U_i$.
By Lemma \ref{lem:flagcomponent}(ii) $\tau_i(A)V$ is equal to
$\sum_{j=i}^d E_jV$, so $\tau_i(A)E^*_0V$ is contained in
$\sum_{j=i}^d E_jV$.
By (\ref{eq:Atrid}) and since  $\tau_i(A)$ has degree $i$
we find $\tau_i(A)E^*_0V$ is contained in $\sum_{k=0}^i E^*_kV$.
By these comments and the definition of $U_i$ we find 
$\tau_i(A)E^*_0V \subseteq U_i$. 
By Lemma \ref{lem:XrEs0Xs} we have $\tau_i(A)E^*_0 \not=0$ so
$\tau_i(A)E^*_0V \not=0$.
By this and since $U_i$ has dimension $1$ we find
$\tau_i(A)E^*_0V=U_i$. We have now shown that
$\tau_i(A)E^*_0V$ is equal to the $i^\text{th}$
component of $[0^*D]$.
In a similar way we find that $\eta^*_{d-i}(A^*)E_dV$ is equal to the
$i^\text{th}$ component of $[0^*D]$.

(ii)--(iv): Apply (i) to the relatives of $\Phi$.
\end{proof}

\section{The action of $S$ on the flags}

Let $(A; \{E_i\}_{i=0}^d; A^*; \{E^*_i\}_{i=0}^d)$ denote 
a Leonard system in $\cal A$
and let $S$ denote the corresponding switching element.
In this section we characterize $S$ via its action on the flags from 
Definition \ref{def:flags}.

\medskip

\begin{theorem}                \label{thm:Sactsflags}
Let $(A; \{E_i\}_{i=0}^d; A^*; \{E^*_i\}_{i=0}^d)$
denote a Leonard system in $\cal A$ and let $S$ denote the corresponding switching element.
Then for all nonzero $X \in {\cal A}$ the following (i), (ii) are equivalent.
\begin{itemize}
\item[(i)] $X$ is a scalar multiple of $S$.
\item[(ii)] $X[0] \subseteq [0]$, $X[D] \subseteq [D]$, 
and  $X[0^*] \subseteq [D^*]$.
\end{itemize}
Suppose (i), (ii) hold. Then equality is attained everywhere in (ii).
\end{theorem}

\begin{proof}
(i)$\Rightarrow$(ii): We show $S[0]=[0]$, $S[D]=[D]$,
and $S[0^*]=[D^*]$.
For $0 \leq i \leq d$ we find
$SE_iV\subseteq E_iV$ since $S$ is a polynomial in $A$,
and $SE_iV=E_iV$ since $S^{-1}$ exists.
Therefore $S[0]=[0]$ and $S[D]=[D]$.
We now show that $S[0^*]=[D^*]$.
To this end we fix an
integer $i$ $(0 \leq i \leq d)$ and show
\begin{equation}        \label{eq:SEs0VEsiV}
  S(E^*_0V+E^*_1V+\cdots+E^*_iV) = E^*_dV+E^*_{d-1}V+\cdots+E^*_{d-i}V.
\end{equation}
Using Lemma \ref{lem:uiAEs0V} and Theorem \ref{thm:SudA} we find that
for $0 \leq j \leq i$,
\begin{eqnarray*}
 SE^*_jV &=& u_d(A)u_j(A)E^*_0V \\
       &=& u_j(A)u_d(A)E^*_0V \\
       &=& u_j(A)E^*_d V. 
\end{eqnarray*}
By (\ref{eq:Atrid}) and since the polynomial $u_j$ has degree $j$,
\[
   u_j(A) E^*_d V \subseteq
     E^*_{d-j}V+E^*_{d-j+1}V+\cdots+E^*_dV.
\]
Combining these comments we obtain
\[
   SE^*_jV \subseteq  E^*_{d-j}V+E^*_{d-j+1}V+\cdots+E^*_dV,
\]
and it follows that
\[
  S(E^*_0V+E^*_1V+\cdots+E^*_iV) \subseteq E^*_dV+E^*_{d-1}V+\cdots+E^*_{d-i}V.
\]
In the above inclusion each side has the same dimension since $S^{-1}$ exists,
so the inclusion becomes equality and (\ref{eq:SEs0VEsiV})
holds.

(ii)$\Rightarrow$(i):
By Theorem \ref{thm:XEs0V} it suffices to
show that $X \in {\cal D}$ and $XE^*_0V \subseteq E^*_dV$.
We first show $X \in {\cal D}$.
Recall that the flags $[0]$, $[D]$ induce the decomposition $[0D]$.
We assume $X[0] \subseteq [0]$ and
$X[D] \subseteq [D]$ so $X [0D] \subseteq [0D]$.
This means that $XE_iV \subseteq E_iV$ for $0 \leq i \leq d$,
so $X \in {\cal D}$.
To get $XE^*_0V \subseteq E^*_dV$, consider the $0^\text{th}$
component in the inclusion $X [0^*] \subseteq [D^*]$.

Suppose (i), (ii) hold. We mentioned in the proof of
(i)$\Rightarrow$(ii) that $S[0]=[0]$, $S[D]=[D]$, and $S[0^*]=[D^*]$.
But $X$ is nonzero and a scalar multiple of $S$ so
$X[0]=[0]$, $X[D]=[D]$, and $X[0^*]=[D^*]$.
\end{proof}

\section{The action of $S$ on the decompositions}

Let $(A; \{E_i\}_{i=0}^d; A^*; \{E^*_i\}_{i=0}^d)$ denote 
a Leonard system in $\cal A$
and let $S$ denote the corresponding switching element.
In this section we characterize
$S$ via its action on the decompositions from
Definition \ref{def:decompositions}.

\medskip

\begin{theorem}      \label{thm:Schar}      \samepage
Let $(A; \{E_i\}_{i=0}^d; A^*; \{E^*_i\}_{i=0}^d)$
denote a Leonard system in $\cal A$ and let $S$ denote the
corresponding switching element.
Then for all nonzero $X \in {\cal A}$
the following (i), (ii) are equivalent.
\begin{itemize}
\item[(i)] $X$ is a scalar multiple of $S$.
\item[(ii)] $X[0^*0] \subseteq [D^*0]$ and $X[0^*D] \subseteq [D^*D]$.
\end{itemize}
Suppose (i), (ii) hold. Then equality holds everywhere in (ii).
\end{theorem}

\begin{proof}
(i)$\Rightarrow$(ii):
By Theorem \ref{thm:Sactsflags}
we have $S[0^*] \subseteq [D^*]$ and $S[0] \subseteq [0]$ so 
$S[0^*0] \subseteq [D^*0]$. Since each component of a decomposition
has dimension $1$ and since $S^{-1}$ exists, we find
$S[0^*0]=[D^*0]$.
In a similar way we obtain $S[0^*D] = [D^*D]$.

(ii)$\Rightarrow$(i):
By Theorem \ref{thm:Sactsflags}
it suffices to show $X[0] \subseteq [0]$, $X[D] \subseteq [D]$,
and $X[0^*] \subseteq [D^*]$.
We first show $X[0] \subseteq [0]$.
By Lemma \ref{lem:induce}(i) and since $X[0^*0] \subseteq [D^*0]$ we find
$X[00^*] \subseteq [0D^*]$.
The decompositions $[00^*]$ and $[0D^*]$ each induce the flag $[0]$ 
by Lemma \ref{lem:induce}(iii) so $X[0] \subseteq [0]$.
Next we show $X[D] \subseteq [D]$.
By Lemma \ref{lem:induce}(i) and since $X[0^*D] \subseteq [D^*D]$ we find
$X[D0^*] \subseteq [DD^*]$.
The decompositions $[D0^*]$ and $[DD^*]$ each induce the flag $[D]$ 
by Lemma \ref{lem:induce}(iii) so $X[D] \subseteq [D]$.
Finally we show $X[0^*] \subseteq [D^*]$.
By Lemma \ref{lem:induce}(iii) we find that 
$[0^*D]$ induces $[0^*]$ and $[D^*D]$ induces $[D^*]$.
By this and since $X[0^*D] \subseteq [D^*D]$ we find
$X[0^*] \subseteq [D^*]$.

Suppose (i), (ii) hold. We mentioned in the proof of
(i)$\Rightarrow$(ii) that $S[0^*0]=[D^*0]$ and $S[0^*D]=[D^*D]$.
But $X$ is nonzero and a scalar multiple of $S$ so
$X[0^*0]=[D^*0]$ and $X[0^*D]=[D^*D]$.
\end{proof}

\section{Some group commutators}

Let $(A; \{E_i\}_{i=0}^d; A^*; \{E^*_i\}_{i=0}^d)$ denote a Leonard system
with switching element $S$ and dual switching element $S^*$.
In this section we consider linear transformations
such as $S^*S^{-1}S^{*-1}S$.
As we will see, these maps are closely
related to the decompositions from Definition \ref{def:decompositions}.
We start with a lemma.

\medskip

\begin{lemma}        \label{lem:fix[0*]}         \samepage
Let $(A; \{E_i\}_{i=0}^d; A^*; \{E^*_i\}_{i=0}^d)$
denote a Leonard system in $\cal A$, with switching element $S$
and dual switching element $S^*$.
Then referring to Definition \ref{def:flags}
the following (i)--(iv) hold.
\begin{itemize}
\item[(i)] $S^*S^{-1}S^{*-1}S$ fixes each of $[0^*]$, $[D]$.
\item[(ii)] $S^*SS^{*-1}S^{-1}$ fixes each of $[D^*]$, $[D]$.
\item[(iii)] $S^{*-1}S^{-1}S^*S$ fixes each of $[0^*]$, $[0]$.
\item[(iv)] $S^{*-1}SS^*S^{-1}$ fixes each of $[D^*]$, $[0]$.
\end{itemize}
\end{lemma}

\begin{proof}
By Theorem \ref{thm:Sactsflags} we find
$S[0]=[0]$, $S[D]=[D]$, and $S[0^*]=[D^*]$.
Applying this to $\Phi^*$ we find
$S^*[0^*]=[0^*]$, $S^*[D^*]=[D^*]$, and $S^*[0]=[D]$.
Combining these comments we routinely obtain the result.
\end{proof}

\medskip

\begin{corollary}       \label{cor:fix[0*d]}            \samepage
Let $(A; \{E_i\}_{i=0}^d; A^*; \{E^*_i\}_{i=0}^d)$ 
denote a Leonard system in $\cal A$, with switching element $S$
and dual switching element $S^*$.
Then referring to Definition \ref{def:decompositions}
the following (i)--(iv) hold.
\begin{itemize}
\item[(i)] $S^*S^{-1}S^{*-1}S$ fixes $[0^*D]$.
\item[(ii)] $S^*SS^{*-1}S^{-1}$ fixes $[D^*D]$.
\item[(iii)] $S^{*-1}S^{-1}S^*S$ fixes $[0^*0]$.
\item[(iv)] $S^{*-1}SS^*S^{-1}$ fixes $[D^*0]$.
\end{itemize}
\end{corollary}

\begin{proof}
(i): For notational convenience abbreviate
$T= S^*S^{-1}S^{*-1}S$.
By Lemma \ref{lem:fix[0*]} we have $T[0^*]=[0^*]$ and $T[D]=[D]$ so
$T[0^*D] \subseteq [0^*D]$.
By this and since $T^{-1}$ exists we find $T[0^*D]=[0^*D]$.

(ii)--(iv) Apply (i) to the relatives of $\Phi$ and
use Theorem \ref{thm:Srelative}.
\end{proof}

\medskip

Referring to Corollary \ref{cor:fix[0*d]}, 
each part (i)--(iv) is asserting that for $0 \leq i \leq d$, 
the $i^\text{th}$ component of the given decomposition is an 
eigenspace for the given operator. 
We now find the corresponding eigenvalue.
We will focus on case (i); the eigenvalues for the remaining cases
will be found using the $D_4$ action.

\medskip

\begin{lemma}    \cite[Theorem 5.2]{NT:mu}   \label{lem:mu1}   \samepage
Let $(A; \{E_i\}_{i=0}^d; A^*; \{E^*_i\}_{i=0}^d)$
denote a Leonard system and let
$(\{\theta_i\}_{i=0}^d; \{\theta^*_i\}_{i=0}^d;
        \{\varphi_i\}_{i=1}^d; \{\phi_i\}_{i=1}^d)$
denote the corresponding parameter array.
Then for $0 \leq i \leq d$,
\begin{eqnarray}
\eta_i(A)E^*_0E_0 &=&
 \frac{\phi_{1}\phi_{2} \cdots \phi_{i}}
      {\eta^*_d(\theta^*_0)} 
         \eta^*_{d-i}(A^*)E_0, 
                                \label{eq:etaiAEs0E0}  \\
\eta_i(A)E^*_dE_0 &=&
 \frac{\varphi_d \varphi_{d-1} \cdots \varphi_{d-i+1}}
      {\tau^*_d (\theta^*_d)}
           \tau^*_{d-i}(A^*)E_0, 
                                 \label{eq:etaiAEsdE0}  \\
\tau_i(A) E^*_0 E_d &=&
 \frac{\varphi_{1}\varphi_{2} \cdots \varphi_{i}}
      {\eta^*_d(\theta^*_0)} 
         \eta^*_{d-i}(A^*)E_d, 
                                 \label{eq:tauiAEs0Ed}   \\
\tau_i(A) E^*_d E_d &=&
 \frac{\phi_{d}\phi_{d-1} \cdots \phi_{d-i+1}}
      {\tau^*_d(\theta^*_d)} 
         \tau^*_{d-i}(A^*)E_d,  
                                   \label{eq:tauiAEsdEd}\\
\eta^*_i(A^*) E_0 E^*_0 &=&
 \frac{\phi_{d}\phi_{d-1} \cdots \phi_{d-i+1}}
      {\eta_d(\theta_0)} 
         \eta_{d-i}(A)E^*_0,
                                  \label{eq:etasiAsE0Es0} \\
\eta^*_i(A^*) E_d E^*_0 &=&
 \frac{\varphi_{d}\varphi_{d-1} \cdots \varphi_{d-i+1}}
      {\tau_d(\theta_d)} 
         \tau_{d-i}(A)E^*_0, 
                                  \label{eq:etasiAsEdEs0} \\
\tau^*_i(A^*) E_0 E^*_d &=&
 \frac{\varphi_{1}\varphi_{2} \cdots \varphi_{i}}
      {\eta_d(\theta_0)} 
         \eta_{d-i}(A)E^*_d, 
                                  \label{eq:tausiAsE0Esd} \\
\tau^*_i(A^*) E_d E^*_d &=&
 \frac{\phi_{1}\phi_{2} \cdots \phi_{i}}
      {\tau_d(\theta_d)} 
         \tau_{d-i}(A)E^*_d.
                                    \label{eq:tausiAsEdEsd}
\end{eqnarray}
\end{lemma}

\medskip

\begin{lemma}   \cite[Theorem 5.6]{NT:mu}    \label{lem:mu2}    \samepage
Let $(A; \{E_i\}_{i=0}^d; A^*; \{E^*_i\}_{i=0}^d)$
denote a Leonard system and let
$(\{\theta_i\}_{i=0}^d; \{\theta^*_i\}_{i=0}^d;
        \{\varphi_i\}_{i=1}^d; \{\phi_i\}_{i=1}^d)$
denote the corresponding parameter array.
Then
\begin{equation}           \label{eq:E0EsdEdEs0}
E_0 E^*_d E_d E^*_0 =
  \frac{\varphi_1 \varphi_2 \cdots \varphi_d}
       {\tau_d(\theta_d) \tau^*_d(\theta^*_d)} E_0 E^*_0.  
\end{equation}
\end{lemma}

\medskip

\begin{lemma}       \label{lem:SEs0}      \samepage
Let $(A; \{E_i\}_{i=0}^d; A^*; \{E^*_i\}_{i=0}^d)$
denote a Leonard system with switching element $S$
and dual switching element $S^*$.
Then the following (i)--(iv) hold.
\begin{itemize}
\item[(i)] 
$SE^*_0$ is equal to each of 
\begin{equation}          \label{eq:SEs0}
  \frac{\tau_d(\theta_d)\tau^*_d(\theta^*_d)}
       {\varphi_1\varphi_2\cdots\varphi_d}   E^*_dE_dE^*_0,
 \qquad\qquad
 \frac{\eta_d(\theta_0)\tau^*_d(\theta^*_d)}
      {\varphi_1\varphi_2\cdots\varphi_d}  E^*_dE_0E^*_0.
\end{equation}
\item[(ii)]
$S^{-1}E^*_d$ is equal to each of
\begin{equation}            \label{eq:SinvEsd} 
  \frac{\tau_d(\theta_d)\eta^*_d(\theta^*_0)}
       {\phi_1\phi_2\cdots\phi_d}  E^*_0E_dE^*_d,
\qquad\qquad
  \frac{\eta_d(\theta_0)\eta^*_d(\theta^*_0)}
       {\phi_1\phi_2\cdots\phi_d}  E^*_0E_0E^*_d.
\end{equation}
\item[(iii)]
$S^* E_0$ is equal to each of
\begin{equation}              \label{eq:SsE0}
  \frac{\tau^*_d(\theta^*_d)\tau_d(\theta_d)}
       {\varphi_1\varphi_2\cdots\varphi_d}  E_dE^*_dE_0,
 \qquad\qquad
  \frac{\eta^*_d(\theta^*_0)\tau_d(\theta_d)}
       {\varphi_1\varphi_2\cdots\varphi_d} E_dE^*_0E_0.
\end{equation}
\item[(iv)]
$S^{*-1}E_d$ is equal to each of
\begin{equation}                 \label{eq:SsinvEd}
  \frac{\tau^*_d(\theta^*_d)\eta_d(\theta_0)}
       {\phi_1\phi_2\cdots\phi_d} E_0E^*_dE_d,
 \qquad\qquad
  \frac{\eta^*_d(\theta^*_0)\eta_d(\theta_0)}
       {\phi_1\phi_2\cdots\phi_d}    E_0E^*_0E_d.
\end{equation}
\end{itemize}
\end{lemma}

\begin{proof}
We first show that $SE^*_0$ is equal to the expression on the left
in (\ref{eq:SEs0}).
By Theorem \ref{thm:XEs0V} we have $SE^*_0V=E^*_dV$ so $SE^*_0=E^*_dSE^*_0$.
The element $E^*_dE_dE^*_0$ is nonzero by (\ref{eq:E0EsdEdEs0})
and Lemma \ref{lem:XrEs0Xs}, 
so it forms a basis for $E^*_d{\cal A}E^*_0$.
This space contains $SE^*_0$ so there exists $\alpha \in \mathbb{K}$
such that $SE^*_0=\alpha E^*_dE_dE^*_0$.
To find $\alpha$, note that $E_0S=E_0$ by (\ref{eq:defS})
so $E_0E^*_0=\alpha E_0E^*_dE_dE^*_0$.
Comparing this with (\ref{eq:E0EsdEdEs0}) we find
\[
   \alpha = \frac{\tau_d(\theta_d)\tau^*_d(\theta^*_d)}
                 {\varphi_1\varphi_2\cdots\varphi_d}.
\]
We have now shown that $SE^*_0$ is equal to the expression
on the left in (\ref{eq:SEs0}).
To obtain the remaining assertions, apply $D_4$
and use Lemma \ref{lem:D4} and Theorem \ref{thm:Srelative}.
\end{proof}

\medskip

\begin{lemma}       \label{lem:prod}   \samepage
Let $(A; \{E_i\}_{i=0}^d; A^*; \{E^*_i\}_{i=0}^d)$
denote a Leonard system with switching element $S$
and dual switching element $S^*$.
Then for $0 \leq i \leq d$,
\begin{equation}             \label{eq:prod} 
S^*S^{-1}S^{*-1}S\tau_i(A)E^*_0 =
 \frac{\phi_1\phi_2\cdots\phi_i}
      {\varphi_1\varphi_2\cdots\varphi_i}
 \frac{\varphi_1\varphi_2\cdots\varphi_{d-i}}
      {\phi_1\phi_2\cdots\phi_{d-i}}
  \tau_i(A)E^*_0.
\end{equation}
\end{lemma}

\begin{proof}
We evaluate the expression on the left in (\ref{eq:prod}).
Recall $S$, $A$ commute; pull $S$ to the right past $\tau_i(A)$.
Now evaluate $SE^*_0$ using the expression on the left in (\ref{eq:SEs0}) 
and in the resulting expression evaluate $\tau_i(A)E^*_dE_d$ 
using (\ref{eq:tauiAEsdEd}); we find the left-hand side of (\ref{eq:prod})
is a scalar multiple of
\begin{equation}           \label{eq:prod1}
  S^*S^{-1}S^{*-1}\tau^*_{d-i}(A^*)E_dE^*_0.
\end{equation}
In line (\ref{eq:prod1}) pull $S^{*-1}$ to the right past $\tau^*_{d-i}(A^*)$.
Now evaluate $S^{*-1}E_d$ using the expression on the left in
(\ref{eq:SsinvEd}) and in the resulting expression
evaluate $\tau^*_{d-i}(A^*)E_0E^*_d$ using (\ref{eq:tausiAsE0Esd});
this shows (\ref{eq:prod1}) is a scalar multiple of 
\begin{equation}           \label{eq:prod2}
  S^*S^{-1}\eta_i(A)E^*_dE_dE^*_0.
\end{equation}
In line (\ref{eq:prod2}) pull $S^{-1}$ to the right past $\eta_i(A)$.
Now evaluate $S^{-1}E^*_d$ using the expression on the right
in (\ref{eq:SinvEsd}) and in the resulting expression
evaluate $\eta_i(A)E^*_0E_0$ using (\ref{eq:etaiAEs0E0});
this shows (\ref{eq:prod2}) is a scalar multiple of
\begin{equation}          \label{eq:prod3}
 S^*\eta^*_{d-i}(A^*)E_0E^*_dE_dE^*_0.
\end{equation}
In line (\ref{eq:prod3}) pull $S^*$ to the right past $\eta^*_{d-i}(A^*)$.
Now evaluate $S^*E_0E^*_dE_d$ using the expression on the left in
(\ref{eq:SsinvEd}) and in the resulting expression
evaluate $\eta^*_{d-i}(A^*)E_dE^*_0$ using (\ref{eq:etasiAsEdEs0});
this shows (\ref{eq:prod3}) is a scalar multiple of $\tau_i(A)E^*_0$.
By the above comments we find that the left-hand side of (\ref{eq:prod})
is a scalar multiple of $\tau_i(A)E^*_0$.
Keeping track of the scalar we routinely verify (\ref{eq:prod}).
\end{proof}

\medskip

\begin{theorem}    \label{thm:main2}   \samepage
Let $(A; \{E_i\}_{i=0}^d; A^*; \{E^*_i\}_{i=0}^d)$
denote a Leonard system in $\cal A$ with parameter array
$(\{\theta_i\}_{i=0}^d; \{\theta^*_i\}_{i=0}^d;
        \{\varphi_i\}_{i=1}^d; \{\phi_i\}_{i=1}^d)$,
switching element $S$ and dual switching element $S^*$.
Then the following (i)--(iv) hold for $0 \leq i \leq d$.
\begin{itemize}
\item[(i)]
The eigenvalue of $S^*S^{-1}S^{*-1}S$ on the $i^\text{th}$ component
of $[0^*D]$ is
\begin{equation}    \label{eq:main}
  \frac{\phi_1\phi_2\cdots\phi_i}
         {\varphi_1\varphi_2\cdots\varphi_i} 
    \frac{\varphi_1\varphi_2\cdots\varphi_{d-i}}
         {\phi_1\phi_2\cdots\phi_{d-i}}.
\end{equation}
\item[(ii)]
The eigenvalue of $S^*SS^{*-1}S^{-1}$ on the $i^\text{th}$ component
of $[D^*D]$ is
\begin{equation}      \label{eq:maindown}
 \frac{\varphi_d\varphi_{d-1}\cdots\varphi_{d-i+1}}
      {\phi_d\phi_{d-1}\cdots\phi_{d-i+1}}
 \frac{\phi_d\phi_{d-1}\cdots\phi_{i+1}}
      {\varphi_d\varphi_{d-1}\cdots\varphi_{i+1}}.
\end{equation}
\item[(iii)]
The eigenvalue of $S^{*-1}S^{-1}S^{*}S$ on the $i^\text{th}$ component
of $[0^*0]$ is
\begin{equation}      \label{eq:mainDown}
 \frac{\varphi_1\varphi_2\cdots\varphi_i}
      {\phi_1\phi_2\cdots\phi_i}
 \frac{\phi_1\phi_2\cdots\phi_{d-i}}
      {\varphi_1\varphi_2\cdots\varphi_{d-i}}.
\end{equation}
\item[(iv)]
The eigenvalue of $S^{*-1}SS^{*}S^{-1}$ on the $i^\text{th}$ component
of $[D^*0]$ is
\begin{equation}      \label{eq:maindownDown}
 \frac{\phi_d\phi_{d-1}\cdots\phi_{d-i+1}}
      {\varphi_d\varphi_{d-1}\cdots\varphi_{d-i+1}}
 \frac{\varphi_d\varphi_{d-1}\cdots\varphi_{i+1}}
      {\phi_d\phi_{d-1}\cdots\phi_{i+1}}.
\end{equation}
\end{itemize}
\end{theorem}

\begin{proof}
(i):
Let $\varepsilon_i$ denote the expression in (\ref{eq:main}).
We show $S^*S^{-1}S^{*-1}S - \varepsilon_i I$ is zero on the $i^\text{th}$
component of $[0^*D]$. But this is immediate from
Lemma \ref{lem:prod} and since this $i^\text{th}$ component equals
$\tau_i(A)E^*_0V$ by Lemma \ref{lem:split}(i).

(ii)--(iv): Apply the $D_4$ action and
use Lemma \ref{lem:D4}, Theorem \ref{thm:Srelative}.
\end{proof}

\section{Representing the elements $S$, $S^*$, $S^{-1}$, $S^{*-1}$ by matrices}

\begin{definition}          \label{def:natural}              \samepage
Let $(A; \{E_i\}_{i=0}^d; A^*; \{E^*_i\}_{i=0}^d)$
denote a Leonard system in $\cal A$ and
fix a nonzero $v^*_0 \in E^*_0V$.
By Lemma \ref{lem:split}(i)
the vectors $\tau_i(A)v^*_0$ $(0 \leq i \leq d)$ form a basis for $V$.
For $X \in {\cal A}$ let $X^{\natural}$ denote the matrix in $\Mat{d+1}$
that represents $X$ with respect to this basis.
We observe $\natural: {\cal A} \to \Mat{d+1}$ is an isomorphism
of $\mathbb{K}$-algebras.
\end{definition}

\medskip

\begin{example}   \cite[Section 21]{T:survey}   \label{exm:AAs}    \samepage
Let $\Phi=(A; \{E_i\}_{i=0}^d; A^*; \{E^*_i\}_{i=0}^d)$
denote a Leonard system in $\cal A$ and
let the isomorphism $\natural : {\cal A} \to \Mat{d+1}$ be as in 
Definition \ref{def:natural}.
Then
\[
A^{\natural} =
\begin{pmatrix}
   \theta_0 &          & & & & \text{\bf 0} \\
   1        & \theta_1 \\
            & 1        & \theta_2 \\
            &          &  \cdot   &  \cdot \\
            &          &          & \cdot  & \cdot\\
   \text{\bf 0} &      &          &       & 1 & \theta_d
\end{pmatrix},
\qquad
A^{*\natural} =
\begin{pmatrix}
   \theta^*_0 &  \varphi_1  & & & & \text{\bf 0} \\
           & \theta^*_1 & \varphi_2 \\
            &         & \theta^*_2 & \cdot \\
            &          &       &  \cdot & \cdot \\
            &          &          &   & \cdot & \varphi_d\\
   \text{\bf 0} &      &          &       &  & \theta^*_d
\end{pmatrix},
\]
where 
$(\{\theta_i\}_{i=0}^d; \{\theta^*_i\}_{i=0}^d;
        \{\varphi_i\}_{i=1}^d; \{\phi_i\}_{i=1}^d)$
denotes the parameter array of $\Phi$.
\end{example}

\medskip

Let $\Phi=(A; \{E_i\}_{i=0}^d; A^*; \{E^*_i\}_{i=0}^d)$
denote a Leonard system in $\cal A$, with switching element $S$ and
dual switching element $S^*$.
Our goal for this section is to find
$S^\natural$, $(S^{-1})^\natural$, $S^{*\natural}$, $(S^{* -1})^\natural$.

\medskip

\begin{lemma}         \label{lem:SetasAsEd}   \samepage
Let $(A; \{E_i\}_{i=0}^d; A^*; \{E^*_i\}_{i=0}^d)$
denote a Leonard system with switching element $S$ and
dual switching element $S^*$.
Then for $0 \leq i \leq d$,
\begin{eqnarray}          
  S\eta^*_{d-i}(A^*)E_d &=&
   \frac{\phi_d\phi_{d-1}\cdots\phi_{d-i+1}}
        {\varphi_1\varphi_2\cdots\varphi_i}
    \tau^*_{d-i}(A^*)E_d,                       \label{eq:SetasEd}  
\\
  S\eta^*_{i}(A^*)E_0 &=&
   \frac{\phi_d\phi_{d-1}\cdots\phi_{d-i+1}}
        {\varphi_1\varphi_2\cdots\varphi_i}
    \tau^*_{i}(A^*)E_0,                       \label{eq:SetasE0} 
\\
  S^* \eta_{d-i}(A)E^*_d &=&
   \frac{\phi_1\phi_{2}\cdots\phi_{i}}
        {\varphi_1\varphi_2\cdots\varphi_i}
    \tau_{d-i}(A)E^*_d,                       \label{eq:SsetaEsd} 
\\
  S^* \eta_{i}(A)E^*_0 &=&
   \frac{\phi_1\phi_{2}\cdots\phi_{i}}
        {\varphi_1\varphi_2\cdots\varphi_i}
    \tau_{i}(A)E^*_0.                       \label{eq:SsetaEs0} 
\end{eqnarray}
\end{lemma}

\begin{proof}
We first show (\ref{eq:SetasEd}).
Using (\ref{eq:tauiAEs0Ed}) we find that
the left-hand side of (\ref{eq:SetasEd}) is a scalar multiple of
\begin{equation}       \label{eq:aaux1}
  S\tau_i(A)E^*_0E_d.
\end{equation}
In (\ref{eq:aaux1}) pull $S$ to the right past $\tau_i(A)$.
Now evaluate $SE^*_0$ using the expression on the left in
(\ref{eq:SEs0}) and in the resulting expression
evaluate $\tau_i(A)E^*_dE_d$ using (\ref{eq:tauiAEsdEd});
this shows that (\ref{eq:aaux1}) is a scalar multiple of
\begin{equation}           \label{eq:aaux2}
  \tau^*_{d-i}(A^*)E_dE^*_0E_d.
\end{equation}
By \cite[Theorem 23.8]{T:survey},
\[
  E_0E^*_0E_0 =
   \frac{\phi_1\phi_2\cdots\phi_d}
        {\eta_d(\theta_0)\eta^*_d(\theta^*_0)} E_0.
\]
Applying $\Downarrow$ to this and using Lemma \ref{lem:D4} we find
\[
  E_dE^*_0E_d = 
   \frac{\varphi_1\varphi_2\cdots\varphi_d}
        {\tau_d(\theta_d)\eta^*_d(\theta^*_0)}E_d.
\]
Using this we find that
(\ref{eq:aaux2}) is a scalar multiple of $\tau^*_{d-i}(A^*)E_d$.
By the above comments the left-hand side of (\ref{eq:SetasEd})
is a scalar multiple of $\tau^*_{d-i}(A^*)E_d$.
Keeping track of the scalar we routinely verify (\ref{eq:SetasEd}).
To obtain (\ref{eq:SetasE0})--(\ref{eq:SsetaEs0}) apply $D_4$ to
(\ref{eq:SetasEd}) using Lemma \ref{lem:D4} and Theorem \ref{thm:Srelative}.
\end{proof}

\medskip

Before we proceed we recall some scalars.
Given a Leonard system 
$(A$; $\{E_i\}_{i=0}^d$; $A^*$; $\{E^*_i\}_{i=0}^d)$
and given nonnegative integers $r,s,t$ such that $r+s+t \leq d$,
in \cite[Definition 13.1]{T:24points} we defined a scalar 
$[r,s,t]_q \in \mathbb{K}$,
where $q+q^{-1}+1$ is the common value of (\ref{eq:indep}).
For example, if $q \neq 1$ and $q \neq -1$ then
\[
  [r,s,t]_q =
   \frac{(q;q)_{r+s}(q;q)_{r+t}(q;q)_{s+t}}
        {(q;q)_r (q;q)_s (q;q)_t (q;q)_{r+s+t}},
\]
where 
\[
     (a;q)_n = (1-a)(1-aq)\cdots(1-aq^{n-1}).
\]
We mention some features of $[r,s,t]_q$ that we will use.
By \cite[Lemma 13.2]{T:24points} we find
$[r,s,t]_q$ is symmetric in $r,s,t$.
We also have the following.

\medskip

\begin{lemma}              \label{lem:X}
Referring to Definition \ref{def:tau}, for  $0 \leq j \leq d$ we have
\begin{eqnarray}
 \tau_j &=& \sum_{i=0}^j [i,j-i,d-j]_q \tau_{j-i}(\theta_d) \eta_i,
               \label{eq:tauj}    \\
 \eta_j &=& \sum_{i=0}^j [i,j-i,d-j]_q \eta_{j-i}(\theta_0) \tau_i,
               \label{eq:etaj}    \\
 \tau^*_j &=& \sum_{i=0}^j [i,j-i,d-j]_q \tau^*_{j-i}(\theta^*_d) \eta^*_i,
               \label{eq:tausj}    \\
 \eta^*_j &=& \sum_{i=0}^j [i,j-i,d-j]_q \eta^*_{j-i}(\theta^*_0) \tau^*_i.
               \label{eq:etasj}
\end{eqnarray}
\end{lemma}

\begin{proof}
This is a routine consequence of \cite[Theorem 15.2]{T:24points}.
\end{proof}

\medskip

\begin{lemma}                 \label{lem:Snatural}            \samepage
Let $(A; \{E_i\}_{i=0}^d; A^*; \{E^*_i\}_{i=0}^d)$
denote a Leonard system with switching element $S$.
Then for $0 \leq j \leq d$,
\begin{eqnarray}
 S \tau_j(A) &=&
  \sum_{i=j}^d
   \frac{[j,i-j,d-i]_q \phi_d\phi_{d-1}\cdots\phi_{d-j+1}\tau^*_{i-j}(\theta^*_d)}
        {\varphi_1\varphi_2\cdots\varphi_i}
       \tau_i(A),
                                     \label{eq:StaujA}   \\
 S^{-1} \tau_j(A) &=&
   \sum_{i=j}^d
     \frac{[j,i-j,d-i]_q \varphi_1\varphi_2\cdots\varphi_j\eta^*_{i-j}(\theta^*_0)}
          {\phi_d\phi_{d-1}\cdots\phi_{d-i+1}}
       \tau_i(A).
                                      \label{eq:SinvtaujA}
\end{eqnarray}
\end{lemma}

\begin{proof}
Concerning (\ref{eq:StaujA}), let $L$ (resp. $R$) denote the expression
on the left (resp. right).
We show $L=R$. To do this we first show that $LE^*_0=RE^*_0$.
We evaluate $LE^*_0$ using (\ref{eq:etasiAsEdEs0}) and
in the resulting expression evaluate 
$S\eta^*_{d-j}(A^*)E_d$ using (\ref{eq:SetasEd});
this shows $LE^*_0$ is a scalar multiple of
\begin{equation}     \label{eq:auxa2}
   \tau^*_{d-j}(A^*)E_dE^*_0.
\end{equation}
Now in (\ref{eq:auxa2}) 
evaluate $\tau^*_{d-j}(A^*)$ using (\ref{eq:tausj}) and
in the resulting expression
evaluate $\eta^*_i(A^*)E_dE^*_0$ using (\ref{eq:etasiAsEdEs0});
this shows that (\ref{eq:auxa2}) is a scalar multiple of
\[
  \sum_{i=0}^{d-j}
    [i,d-j-i,j]_q \varphi_d\varphi_{d-1}\cdots\varphi_{d-i+1}
    \tau^*_{d-j-i}(\theta^*_d) \tau_{d-i}(A) E^*_0.
\]
In this expression we replace $i$ by $d-i$ and find it is equal to
\begin{equation}                \label{eq:aux}
  \sum_{i=j}^d
   [j,i-j,d-i]_q \varphi_d\varphi_{d-1}\cdots\varphi_{i+1}
   \tau^*_{i-j}(\theta^*_d) \tau_i(A)E^*_0.
\end{equation}
So far we have shown that $LE^*_0$ is a scalar multiple of (\ref{eq:aux}).
Keeping track of the scalar we routinely find $LE^*_0=RE^*_0$.
Now $L=R$ by Lemma \ref{lem:XrEs0Xs}.
To obtain (\ref{eq:SinvtaujA}), apply $\downarrow$
to (\ref{eq:StaujA}) and recall $S^{\downarrow}=S^{-1}$ from
Theorem \ref{thm:Srelative}.
\end{proof}

\medskip

\begin{theorem}               \label{thm:Snatural}   \samepage
Let $\Phi=(A; \{E_i\}_{i=0}^d; A^*; \{E^*_i\}_{i=0}^d)$
denote a Leonard system in $\cal A$ and let
$(\{\theta_i\}_{i=0}^d; \{\theta^*_i\}_{i=0}^d;
        \{\varphi_i\}_{i=1}^d; \{\phi_i\}_{i=1}^d)$
denote the corresponding parameter array.
Let $S$ denote the switching element for $\Phi$ and let
the isomorphism $\natural:{\cal A} \to \Mat{d+1}$ be from
Definition \ref{def:natural}.
Then each of $S^{\natural}$, $(S^{-1})^\natural$ is lower triangular.
Moreover for $0 \leq j \leq i \leq d$ their $(i,j)$ entries are given as follows.
\begin{eqnarray*}
S^\natural_{i,j} &=& 
  \frac{[j,i-j,d-i]_q\phi_d\phi_{d-1}\cdots\phi_{d-j+1} \tau^*_{i-j}(\theta^*_d)}
       {\varphi_1\varphi_2\cdots\varphi_i},            \\
(S^{-1})^{\natural}_{i,j} &=& 
  \frac{[j,i-j,d-i]_q \varphi_1\varphi_2\cdots\varphi_j \eta^*_{i-j}(\theta^*_0)}
       {\phi_d\phi_{d-1}\cdots\phi_{d-i+1}}.     
\end{eqnarray*}
\end{theorem}

\begin{proof}
Follows from Lemma \ref{lem:Snatural} and Definition \ref{def:natural}.
\end{proof}

\medskip

\begin{lemma}                 \label{lem:Ssnatural}        \samepage
Let $(A; \{E_i\}_{i=0}^d; A^*; \{E^*_i\}_{i=0}^d)$
denote a Leonard system with dual switching element $S^*$.
Then for $0 \leq j \leq d$,
\begin{eqnarray}
 S^*\tau_j(A)E^*_0 &=&
  \sum_{i=0}^j
    \frac{[i,j-i,d-j]_q \phi_1\phi_2\cdots\phi_i \tau_{j-i}(\theta_d)}
         {\varphi_1\varphi_2\cdots\varphi_i}
      \tau_i(A)E^*_0,                       
                                    \label{eq:SstaujAEs0}    \\
 S^{*-1} \tau_j(A)E^*_0 &=&
  \sum_{i=0}^j
    \frac{[i,j-i,d-j]_q \varphi_1\varphi_2\cdots\varphi_j\eta_{j-i}(\theta_0)}
         {\phi_1\phi_2\cdots\phi_j}
      \tau_i(A)E^*_0.
                                    \label{eq:SsinvtaujAEs0}
\end{eqnarray}
\end{lemma}

\begin{proof}
First we show (\ref{eq:SstaujAEs0}).
In the left-hand side of (\ref{eq:SstaujAEs0}) we evaluate
$\tau_j(A)$ using (\ref{eq:tauj}) to find 
\[
  S^*\tau_j(A)E^*_0 =
   \sum_{i=0}^j [i,j-i,d-j]_q \tau_{j-i}(\theta_d)S^*\eta_i(A)E^*_0.
\]
In this equation we evaluate $S^*\eta_i(A)E^*_0$ using
(\ref{eq:SsetaEs0}) and get (\ref{eq:SstaujAEs0}).

Next we show (\ref{eq:SsinvtaujAEs0}).
By (\ref{eq:SsetaEs0}),
\[
  S^{*-1}\tau_j(A)E^*_0 =
    \frac{\varphi_1\varphi_2\cdots\varphi_j}
         {\phi_1\phi_2\cdots\phi_j}
     \eta_j(A)E^*_0.
\]
In this equation we evaluate $\eta_j(A)$ using (\ref{eq:etaj})
and get (\ref{eq:SsinvtaujAEs0}).
\end{proof}

\medskip

\begin{theorem}               \label{thm:Ssnatural}   \samepage
Let $\Phi=(A; \{E_i\}_{i=0}^d; A^*; \{E^*_i\}_{i=0}^d)$
denote a Leonard system in $\cal A$ and let
$(\{\theta_i\}_{i=0}^d; \{\theta^*_i\}_{i=0}^d;
        \{\varphi_i\}_{i=1}^d; \{\phi_i\}_{i=1}^d)$
denote the corresponding parameter array.
Let $S^*$ denote the dual switching element for $\Phi$ and let
the isomorphism $\natural:{\cal A} \to \Mat{d+1}$ be from
Definition \ref{def:natural}.
Then  each of $S^{*\natural}$, $(S^{*-1})^{\natural}$ is
upper triangular. Moreover for $0 \leq i \leq j \leq d$
their $(i,j)$ entries are given as follows.
\begin{eqnarray*}
 S^{*\natural}_{i,j} &=& 
  \frac{[i,j-i,d-j]_q \phi_1\phi_2\cdots\phi_i \tau_{j-i}(\theta_d)}
       {\varphi_1\varphi_2\cdots\varphi_i},       \\
(S^{*-1})^{\natural}_{i,j} &=& 
  \frac{[i,j-i,d-j]_q \varphi_1\varphi_2\cdots\varphi_j \eta_{j-i}(\theta_0)}
       {\phi_1\phi_2\cdots\phi_j}.     
\end{eqnarray*}
\end{theorem}

\begin{proof}
Follows from Lemma \ref{lem:Ssnatural} and Definition \ref{def:natural}.
\end{proof}

\section{Leonard pairs in matrix form}

In this section we restate Theorem \ref{thm:Snatural} and 
Theorem \ref{thm:Ssnatural} in more concrete terms. 
Let us consider the following situation.

\medskip

\begin{definition}       \label{def:matrices}            \samepage
Let $d$ denote a nonnegative integer.
Let $A$ and $A^*$ denote matrices in $\Mat{d+1}$ of the form
\begin{equation}            \label{eq:AAs2}
A=
 \begin{pmatrix}
   \theta_0 & & & & & \text{\bf 0} \\
   1 & \theta_1 \\
     & 1 & \theta_2 \\
     &   & \cdot & \cdot \\
     &   &       & \cdot & \cdot \\
   \text{\bf 0} & & & & 1 & \theta_d
 \end{pmatrix},
\qquad
A^*=
 \begin{pmatrix}
   \theta^*_0 & \varphi_1 & & & & \text{\bf 0} \\
    & \theta^*_1 & \varphi_2 \\
     &  & \theta^*_2 & \cdot \\
     &   &    & \cdot & \cdot \\
     &   &    &       & \cdot & \varphi_d \\
   \text{\bf 0} & & &  & & \theta^*_d
 \end{pmatrix},
\end{equation}
where
\[
 \theta_i \neq \theta_j, \qquad \theta^*_i \neq \theta^*_j
  \qquad\qquad (0 \leq i,j \leq d),
\]
\[
  \varphi_i \neq 0 \qquad\qquad (1 \leq i \leq d).
\]
Observe $A$ (resp. $A^*$) is multiplicity-free, 
with eigenvalues $\{\theta_i\}_{i=0}^d$ 
(resp. $\{\theta^*_i\}_{i=0}^d$).
For $0 \leq i \leq d$ let $E_i$ (resp. $E^*_i$) denote the
primitive idempotent of $A$ (resp. $A^*$) associated with 
$\theta_i$ (resp. $\theta^*_i$).
\end{definition}

\medskip

\begin{lemma}     \label{lem:new}    \samepage
Referring to Definition \ref{def:matrices}, the following (i), (ii)
are equivalent.
\begin{itemize}
\item[(i)] The pair $A,A^*$ is a Leonard pair.
\item[(ii)] The sequence $\Phi=(A; \{E_i\}_{i=0}^d; A^*; \{E^*_i\}_{i=0}^d)$ 
is a Leonard system.
\end{itemize}
Suppose (i), (ii) hold. Then the sequence $\{\varphi_i\}_{i=1}^d$
is the first split sequence of $\Phi$.
Moreover the isomorphism $\natural$ from Definition \ref{def:natural}
is the identity map.
\end{lemma}

\begin{proof}
The equivalence of (i), (ii) is established in \cite[Lemma 6.2]{T:split}.
By \cite[Theorem 14.3]{T:canform}, 
the sequence $\{\varphi_i\}_{i=1}^d$ is the first split sequence of $\Phi$.
Comparing Example \ref{exm:AAs} with (\ref{eq:AAs2}) 
we find $A^\natural = A$ and $A^{*\natural} = A^*$. 
Now $\natural$ is the identity map
since $A,A^*$ generate $\Mat{d+1}$ by \cite[Corollary 5.5]{T:qRacah}.
\end{proof}

\medskip

\begin{theorem}            \samepage
Referring to Lemma \ref{lem:new},
assume the equivalent conditions
(i), (ii) hold, and let $S$ denote the switching
element element for $\Phi$.
Then each of $S$, $S^{-1}$ is lower triangular.
Moreover for $0 \leq j \leq i \leq d$ their $(i,j)$ entries are given as follows.
\begin{eqnarray*}
S_{i,j} &=& 
  \frac{[j,i-j,d-i]_q\phi_d\phi_{d-1}\cdots\phi_{d-j+1} \tau^*_{i-j}(\theta^*_d)}
       {\varphi_1\varphi_2\cdots\varphi_i},           \\
S^{-1}_{i,j} &=& 
  \frac{[j,i-j,d-i]_q \varphi_1\varphi_2\cdots\varphi_j \eta^*_{i-j}(\theta^*_0)}
       {\phi_d\phi_{d-1}\cdots\phi_{d-i+1}}.    
\end{eqnarray*}
In the above lines $\{\phi_i\}_{i=1}^d$ is the second split sequence of $\Phi$.
\end{theorem}

\begin{proof}
Combine Theorem \ref{thm:Snatural} with Lemma \ref{lem:new}.
\end{proof}

\medskip

\begin{theorem}               \samepage
Referring to Lemma \ref{lem:new},
assume the equivalent conditions
(i), (ii) hold, and let $S^*$ denote the dual switching
element element for $\Phi$.
Then  each of $S^*$, $S^{*-1}$ is
upper triangular. Moreover for $0 \leq i \leq j \leq d$
their $(i,j)$ entries are given as follows.
\begin{eqnarray*}
S^*_{i,j} &=& 
  \frac{[i,j-i,d-j]_q \phi_1\phi_2\cdots\phi_i \tau_{j-i}(\theta_d)}
       {\varphi_1\varphi_2\cdots\varphi_i},        \\
S^{*-1}_{i,j} &=& 
  \frac{[i,j-i,d-j]_q \varphi_1\varphi_2\cdots\varphi_j \eta_{j-i}(\theta_0)}
       {\phi_1\phi_2\cdots\phi_j}.       
\end{eqnarray*}
In the above lines $\{\phi_i\}_{i=1}^d$ is the second split sequence of $\Phi$.
\end{theorem}

\begin{proof}
Combine Theorem \ref{thm:Ssnatural} with Lemma \ref{lem:new}.
\end{proof}

\section{A characterization of a Leonard system in terms of the switching element}

In this section we give a characterization of a Leonard
system in terms of its switching element.
This characterization
is a variation on \cite[Theorem 6.3]{T:split} and is stated as follows.

\medskip

\begin{theorem}          \label{thm:CharLS1}            \samepage
Referring to Definition \ref{def:matrices},
the following (i), (ii) are equivalent.
\begin{itemize}
\item[(i)] The pair $A,A^*$ is a Leonard pair.
\item[(ii)] There exists an invertible $X \in \Mat{d+1}$ and
there exist nonzero scalars $\phi_i \in \mathbb{K}$ $(1 \leq i \leq d)$
such that $X^{-1}AX=A$ and
\begin{equation}         \label{eq:SinvAsS}
X^{-1}A^*X =
 \begin{pmatrix}
   \theta^*_d & \phi_d & & & & \text{\bf 0} \\
    & \theta^*_{d-1} & \phi_{d-1} \\
     &  & \theta^*_{d-2} & \cdot \\
     &   &    & \cdot & \cdot \\
     &   &    &       & \cdot & \phi_1 \\
   \text{\bf 0} & & &  & & \theta^*_0
 \end{pmatrix}.
\end{equation}
\end{itemize}
Suppose (i), (ii) hold. Then
$\Phi=(A; \{E_i\}_{i=0}^d; A^*; \{E^*_i\}_{i=0}^d)$
is a Leonard system, 
$\{\phi_i\}_{i=1}^d$ is the second split sequence of $\Phi$, 
and $X$ is a scalar multiple of the switching element for $\Phi$.
\end{theorem}

\begin{proof}
For notational convenience we abbreviate $V=\mathbb{K}^{d+1}$.
For $0 \leq i \leq d$ let $e_i$ denote the vector in
$V$ with $i^\text{th}$ coordinate $1$ and all other coordinates $0$.
Observe that $\{e_i\}_{i=0}^d$ is a basis for $V$. 
From the form of A in (\ref{eq:AAs2}) we find
\begin{equation}          \label{eq:Aei}
    (A-\theta_iI)e_i =e_{i+1} \quad (0 \leq i \leq d-1),
   \qquad 
    (A-\theta_dI)e_d=0.
\end{equation}
From the form of $A^*$ in (\ref{eq:AAs2}) we find
\begin{equation}        \label{eq:Asei}
   (A^*-\theta^*_iI)e_i = \varphi_i e_{i-1}  \quad (1 \leq i \leq d),
 \qquad
 (A^*-\theta^*_0I)e_0 = 0.
\end{equation}

(i)$\Rightarrow$(ii):
By Lemma \ref{lem:new} the sequence 
$\Phi=(A; \{E_i\}_{i=0}^d; A^*; \{E^*_i\}_{i=0}^d)$
is a Leonard system in $\Mat{d+1}$, with eigenvalue
sequence $\{\theta_i\}_{i=0}^d$, dual eigenvalue sequence $\{\theta^*_i\}_{i=0}^d$,
and first split sequence $\{\varphi_i\}_{i=1}^d$.
Let $\{\phi_i\}_{i=1}^d$ denote the second split sequence
of $\Phi$ and let $S$ denote the switching element for $\Phi$.
Then $S^{-1}$ exists by Lemma  \ref{note:Sinv}.
Also $S$ commutes with $A$ by Theorem \ref{thm:XEs0V}
so $S^{-1}AS=A$.
We now show that (\ref{eq:SinvAsS}) holds with $X=S$.
From (\ref{eq:Aei}) we find
\begin{equation}          \label{eq:ei}
  e_i = \tau_i(A)e_0  \qquad  (0 \leq i \leq d).
\end{equation}
From the equation on the right in (\ref{eq:Asei}) we find
$A^*e_0=\theta^*_0 e_0$; using this and $e_0 \neq 0$ we find
$e_0$ is a basis for $E^*_0V$.
By this and (\ref{eq:Asei}), (\ref{eq:ei}) we have
\begin{equation}           \label{eq:AstauAEs0}
  (A^*-\theta^*_iI)\tau_i(A)E^*_0 = \varphi_i \tau_{i-1}(A)E^*_0 
    \quad(1 \leq i \leq d),
 \qquad
  (A^*-\theta^*_0I)E^*_0 = 0.
\end{equation}
We apply $\downarrow$ to (\ref{eq:AstauAEs0}) and use 
Lemma \ref{lem:D4} to find
\begin{equation}       \label{eq:AstauiAEsd}
  (A^*-\theta^*_{d-i}I)\tau_i(A)E^*_d = 
          \phi_{d-i+1} \tau_{i-1}(A)E^*_d \quad(1 \leq i \leq d),
 \qquad
  (A^*-\theta^*_dI)E^*_d = 0.
\end{equation}
By Theorem \ref{thm:XEs0V} and since $e_0 \in E^*_0V$
we find $Se_0 \in E^*_dV$. Combining this with (\ref{eq:AstauiAEsd})
we have
\begin{equation}        \label{eq:AstauASe0}
  (A^*-\theta^*_{d-i}I)\tau_i(A)Se_0 = 
          \phi_{d-i+1} \tau_{i-1}(A)Se_0 \quad(1 \leq i \leq d),
 \qquad
  (A^*-\theta^*_dI)Se_0 = 0.
\end{equation}
Evaluating (\ref{eq:AstauASe0}) using $S^{-1}AS=A$ and (\ref{eq:ei}) we routinely
find
\[
(S^{-1}A^*S-\theta^*_{d-i}I)e_i = \phi_{d-i+1}e_{i-1}
    \quad  (1 \leq i \leq d),
    \qquad  (S^{-1}A^*S-\theta^*_{d}I)e_0=0.
\]
By this we find (\ref{eq:SinvAsS}) holds with $X=S$, as desired.

(ii)$\Rightarrow$(i):
We show 
$\Phi=(A; \{E_i\}_{i=0}^d; A^*; \{E^*_i\}_{i=0}^d)$
is a Leonard system in $\Mat{d+1}$. 
To do this we invoke \cite[Theorem 5.1]{T:split}.
According to that theorem it suffices to display
a decomposition $\{U_i\}_{i=0}^d$ of $V$ such that
\begin{equation}       \label{eq:Asplit}
 (A-\theta_iI)U_i = U_{i+1} \quad (0 \leq i \leq d-1),
 \qquad
  (A-\theta_dI)U_d=0,
\end{equation}
\begin{equation}        \label{eq:Assplit}
 (A^*-\theta^*_iI)U_i = U_{i-1} \quad   (1 \leq i \leq d),
   \qquad
   (A^*-\theta^*_0I)U_0=0,
\end{equation}
and a decomposition $\{V_i\}_{i=0}^d$ of $V$ such that
\begin{equation}        \label{eq:Asplit2}
  (A-\theta_iI)V_i = V_{i+1}  \quad   (0 \leq i \leq d-1),
  \qquad
   (A-\theta_dI)V_d=0, 
\end{equation}
\begin{equation}       \label{eq:Assplit2}
 (A^*-\theta^*_{d-i}I)V_i = V_{i-1} \quad (1 \leq i \leq d),
   \qquad
   (A^*-\theta^*_dI)V_0=0.
\end{equation}
Define $U_i= \text{Span}\{e_i\}$ for $0 \leq i \leq d$.
The sequence $\{U_i\}_{i=0}^d$ is a decomposition of $V$ since
$\{e_i\}_{i=0}^d$ is a basis for $V$.
The decomposition $\{U_i\}_{i=0}^d$ satisfies
(\ref{eq:Asplit}) by (\ref{eq:Aei}). 
The decomposition $\{U_i\}_{i=0}^d$ satisfies
(\ref{eq:Assplit}) by (\ref{eq:Asei}) and since 
$\varphi_i \neq 0$ for $1 \leq i \leq d$.
Now define $V_i= \text{Span}\{Xe_i\}$ for $0 \leq i \leq d$.
The sequence $\{V_i\}_{i=0}^d$ is a decomposition of $V$ since
$X^{-1}$ exists, and since $\{e_i\}_{i=0}^d$ is a basis for $V$.
The decomposition $\{V_i\}_{i=0}^d$ satisfies (\ref{eq:Asplit2}) by
(\ref{eq:Aei}) and since $X^{-1}AX=A$.
The decomposition $\{V_i\}_{i=0}^d$ satisfies (\ref{eq:Assplit2}) by
(\ref{eq:SinvAsS}) and since $\phi_i \neq 0$ for $1 \leq i \leq d$.
We have now verified (\ref{eq:Asplit})--(\ref{eq:Assplit2}) so 
\cite[Theorem 5.1]{T:split}
applies; by that theorem $\Phi$ is a Leonard system in $\Mat{d+1}$.
Now the pair $A,A^*$ is a Leonard pair by Lemma \ref{lem:new}.

Now suppose (i), (ii) hold. We mentioned in the proof
of (i)$\Rightarrow$(ii) that 
$\Phi=(A; \{E_i\}_{i=0}^d$; $A^*; \{E^*_i\}_{i=0}^d)$
is a Leonard system in $\Mat{d+1}$. Next we show
that $X$ is a scalar multiple of the switching element $S$ for $\Phi$. 
To do this we invoke Theorem \ref{thm:XEs0V}.
Let $\cal D$ denote the subalgebra of $\Mat{d+1}$ generated
by $A$. 
The element $X$ commutes with $A$ so $X \in {\cal D}$.
By the left-most column in (\ref{eq:SinvAsS}) we find
$X^{-1}A^*X e_0= \theta^*_d e_0$ so $Xe_0 \in E^*_dV$.
But $e_0$ is a basis for
$E^*_0V$ so $XE^*_0V \subseteq E^*_dV$.
Now $X$ is a scalar multiple of $S$ by Theorem \ref{thm:XEs0V}.
We saw in the proof of (i)$\Rightarrow$(ii) that the
sequence $\{\phi_i\}_{i=1}^d$ is the second split sequence of $\Phi$.
\end{proof}

\medskip

The following is a variation on \cite[Theorem 3.2]{T:array}.

\medskip

\begin{theorem}       \label{thm:charLS2}   \samepage
Let $d$ denote a nonnegative integer and let
\begin{equation}             \label{eq:parray}
  (\{\theta_i\}_{i=0}^d; \{\theta^*_i\}_{i=0}^d;
        \{\varphi_i\}_{i=1}^d; \{\phi_i\}_{i=1}^d)
\end{equation}
denote a sequence of scalars taken from $\mathbb{K}$.
Assume this sequence satisfies the conditions
(PA1) and (PA2) in Theorem \ref{thm:classify}.
Then the following (i), (ii) are equivalent.
\begin{itemize}
\item[(i)]
The sequence (\ref{eq:parray}) satisfies (PA3)--(PA5) in 
Theorem \ref{thm:classify}.
\item[(ii)]
There exists an invertible $X \in \Mat{d+1}$ such that
\[   
 X^{-1}
 \begin{pmatrix}
   \theta_0 & & & & & \text{\bf 0} \\
   1 & \theta_1 \\
     & 1 & \theta_2 \\
     &   & \cdot & \cdot \\
     &   &       & \cdot & \cdot \\
   \text{\bf 0} & & & & 1 & \theta_d
 \end{pmatrix}
 X =
  \begin{pmatrix}
   \theta_0 & & & & & \text{\bf 0} \\
   1 & \theta_1 \\
     & 1 & \theta_2 \\
     &   & \cdot & \cdot \\
     &   &       & \cdot & \cdot \\
   \text{\bf 0} & & & & 1 & \theta_d
 \end{pmatrix},
\]
\[
 X^{-1}
 \begin{pmatrix}
   \theta^*_0 & \varphi_1 & & & & \text{\bf 0} \\
    & \theta^*_1 & \varphi_2 \\
     &  & \theta^*_2 & \cdot \\
     &   &    & \cdot & \cdot \\
     &   &    &       & \cdot & \varphi_d \\
   \text{\bf 0} & & &  & & \theta^*_d
 \end{pmatrix}
 X =
 \begin{pmatrix}
   \theta^*_d & \phi_d & & & & \text{\bf 0} \\
    & \theta^*_{d-1} & \phi_{d-1} \\
     &  & \theta^*_{d-2} & \cdot \\
     &   &    & \cdot & \cdot \\
     &   &    &       & \cdot & \phi_1 \\
   \text{\bf 0} & & &  & & \theta^*_0
 \end{pmatrix}
\]
\end{itemize}
\end{theorem}

\begin{proof}
(i)$\Rightarrow$(ii):
By Theorem \ref{thm:classify} there exists a Leonard system
$\Phi=(A; \{E_i\}_{i=0}^d; A^*$; $\{E^*_i\}_{i=0}^d)$ over $\mathbb{K}$
that has parameter array (\ref{eq:parray}).
Let $\natural$ denote the corresponding isomorphism from
Definition \ref{def:natural}.
Applying $\natural$ to each term in $\Phi$ if necessary, 
we may assume $\Phi$ is in $\Mat{d+1}$, and
that $\natural$ is the identity map. Now $A,A^*$ are of the form (\ref{eq:AAs2}).
Now (ii) holds by Theorem \ref{thm:CharLS1}.

(ii)$\Rightarrow$(i):
Follows from Theorem \ref{thm:CharLS1} and Theorem \ref{thm:classify}.
\end{proof}

\medskip

\begin{note}
We comment on how the switching element is related to
the matrix $G$ that appears
in \cite[Theorem 6.3]{T:split} and
\cite[Theorem 3.2]{T:array}.
In \cite[Theorem 6.3]{T:split} reference is made to a Leonard pair
$A,A^*$ of the form (\ref{eq:AAs2}); 
let $\Phi=(A; \{E_i\}_{i=0}^d; A^*$; $\{E^*_i\}_{i=0}^d)$ denote
the corresponding Leonard system from Lemma \ref{lem:new}. 
Then $G$ is a nonzero scalar multiple of $S^{*-1}Y$, where
$S^*$ denotes the dual switching element for $\Phi$ and
$Y$ denotes the diagonal matrix in $\Mat{d+1}$ whose $(i,i)$ entry
is
\[
  \frac{\phi_1\phi_2\cdots\phi_i}
       {\varphi_1\varphi_2\cdots\varphi_i}
\]
for $0 \leq i \leq d$. 
\end{note}

\begin{proof}
From the data given in \cite[Theorem 6.3]{T:split} 
one readily verifies that $YG^{-1}$ is invertible and
commutes with $A^*$. 
Moreover $YG^{-1}E_0V \subseteq E_dV$ where $V=\mathbb{K}^{d+1}$.
Now applying Theorem \ref{thm:XEs0V} to $\Phi^*$ we find
$YG^{-1}$ is a nonzero scalar multiple of $S^*$.
So $G$ is a nonzero scalar multiple of $S^{*-1}Y$.
\end{proof}

\section{Open problems}

In this section 
$\Phi=(A; \{E_i\}_{i=0}^d; A^*; \{E^*_i\}_{i=0}^d)$
denotes a Leonard system in $\cal A$ with switching element $S$ 
and dual switching element $S^*$.
We will be discussing the flags
\begin{equation}          \label{eq:z}
  [z]  \qquad (z \in \Omega)
\end{equation}
from Definition \ref{def:flags}.

\medskip

\begin{definition}                \label{def:switch}
For distinct $x,y \in \Omega$, by a 
{\em switching element of type $(x,y)$} we mean 
an element $X$ in ${\cal A}$ that
sends $[x]$ to $[y]$ and fixes the remaining two flags in (\ref{eq:z}).
\end{definition}

\medskip

\begin{example}
By Theorem \ref{thm:Sactsflags} we find that for
nonzero $X \in {\cal A}$,
\begin{itemize}
\item[(i)] 
$X$ is a switching element of type $(0^*,D^*)$ if and only if
$X$ is a scalar multiple of $S$.
\item[(ii)] 
$X$ is a switching element of type $(D^*,0^*)$ if and only if $X$ 
is a scalar multiple of $S^{-1}$.
\item[(iii)]
$X$ is a switching element of type $(0,D)$ if and only if 
$X$ is a scalar multiple of $S^*$.
\item[(iv)]
$X$ is a switching element of type $(D,0)$ if and only if 
$X$ is a scalar multiple of $S^{*-1}$.
\end{itemize}
\end{example}

\medskip

\begin{problem}              \label{prob:1}
Find a necessary and sufficient condition
on the parameter array of $\Phi$
for there to exist a switching element of type $(0^*,D)$.
\end{problem}

\medskip

\begin{problem}              \label{prob:2}
Find a necessary and sufficient condition
on the parameter array of $\Phi$ for there to exist a switching 
element of type $(0^*,D)$ and a switching element of
type $(0,D^*)$.
\end{problem}

\medskip

\begin{note}
For certain $\Phi$ there exists a second Leonard system
$\Phi'= (B; \{F_i\}_{i=0}^d; B^*; \{F^*_i\}_{i=0}^d)$
in $\cal A$ such that the decomposition $\{F_iV\}_{i=0}^d$ 
coincides with $[0D^*]$ and the decomposition 
$\{F^*_iV\}_{i=0}^d$ coincides with $[0^*D]$. 
See for example \cite{H}. 
In this case the switching element for $\Phi'$ is a switching element
for $\Phi$ of type $(0^*,D)$, and the dual switching element for
$\Phi'$ is a switching element for $\Phi$ of type $(0,D^*)$.
\end{note}

\bigskip

\bibliographystyle{plain}

\bigskip\bigskip\noindent
Kazumasa Nomura\\
College of Liberal Arts and Sciences\\
Tokyo Medical and Dental University\\
Kohnodai, Ichikawa, 272-0827 Japan\\
email: nomura.las@tmd.ac.jp

\bigskip\noindent
Paul Terwilliger\\
Department of Mathematics\\
University of Wisconsin\\
480 Lincoln Drive\\ 
Madison, Wisconsin, 53706 USA\\
email: terwilli@math.wisc.edu

\bigskip\noindent
{\bf Keywords.}
Leonard pair, tridiagonal pair, $q$-Racah polynomial, orthogonal polynomial.

\noindent
{\bf 2000 Mathematics Subject Classification}.
05E35, 05E30, 33C45, 33D45.

\end{document}